\documentclass[showpacs,nofootinbib,showkeys,eqsecnum,prd,aps,preprint]{revtex4-1}

\usepackage[english]{babel}
\usepackage[T1]{fontenc}
\usepackage{amssymb}
\usepackage{amsmath}
\usepackage{amsfonts}
\usepackage{latexsym}
\usepackage{bm}
\usepackage{tensor}
\usepackage{hyperref}

\begin{document}

\title{An application of the Maslov complex germ method to the 1D nonlocal Fisher--KPP equation}

\author{A.~V.~Shapovalov${}^{a,b,c}$}
\email{shpv@phys.tsu.ru}

\author{A.~Yu.~Trifonov${}^{b,c}$}
\email{atrifonov@tpu.ru}

\date{\today }

\begin{abstract}

A semiclassical approximation approach based on the Maslov complex germ method
is considered in detail for the 1D  nonlocal Fisher--Kolmogorov--Petrovskii--Piskunov
equation  under the supposition of  weak diffusion. In terms of the semiclassical formalism developed,  the original nonlinear equation is reduced to an associated linear partial differential equation and some algebraic equations for the coefficients of the linear equation with a given accuracy of the asymptotic parameter. The solutions of the nonlinear equation are constructed from the solutions of both the linear equation and the
algebraic equations. The solutions of the linear problem are found with the use of symmetry operators.
A countable family of the leading terms of the semiclassical asymptotics is constructed in explicit form.

The semiclassical asymptotics are valid by construction in a finite time interval. We construct asymptotics  which are different from the semiclassical ones and can describe evolution of the solutions of the Fisher--Kolmogorov--Petrovskii--Piskunov equation at large times.  In the example considered, an initial unimodal  distribution becomes multimodal, which can be treated as  an example of a space structure.

\end{abstract}


\affiliation{
$^{a}$Department of Theoretical Physics, Tomsk State University, Novosobornaya Sq. 1, Tomsk, 634050, Russia\\
${}^{b}$Department of Higher Mathematics and Mathematical Physics, \\
Tomsk Polytechnic University, Lenin Ave., 30, Tomsk, 634034, Russia \\
${}^{c}$Tomsk State  Pedagogical University,\\  Kievskaya St. 60, Tomsk, 634041, Russia
}
\maketitle

\section{Introduction}

Integro-differential equations (IDEs) of reaction--diffusion (RD) type are the theoretical basis for a study of nonlinear phenomena with long-range interactions in condensed matter physics, chemical kinetics, population biology, and other fields of nonlinear science.

Nonlocal versions of the classical Fisher \cite{fisher37} and Kolmogorov--Petrovskii--Piskunov \cite{kpp37} (Fisher--KPP) population equation can describe the nonlinear phenomena inherent in RD systems even for one-species population dynamics, such as formation of population patterns \cite{murray1} and spatially localized lumps \cite{lopez2}, clustering \cite{lopez1}, propagation of traveling and spiral waves \cite{murray1,okubo}, and stability of steady states \cite{kuehn}.

In study of nonlinear phenomena,  pattern formation processes are of a special interest. Spatially or spatiotemporal structures can arise in nonlinear systems being  far from equilibrium when a system is in an unstable state.
In conditions where instability can break homogeneity the system passes into a spatially inhomogeneous stable state which is recognized as a space pattern. Knowledge of these features of nonlinear systems renders better understanding the nature of their evolution. The above examples are essentially particular manifestations of the general pattern formation phenomena.

Mathematical modeling of these phenomena contributes to comprehension of the cancer morphogenesis
(see, e.g., \cite{murray1,murray2,meinhardt}) and the tumor growth \cite{swanson,wang,perez}.
The population dynamics under external influences is of particular interest from the standpoint of possible process control.
The external factors are modeled by the variable coefficients in the integro-differential Fisher--KPP equation.
Equations of this type are mathematically complex, and only few exact integration methods are available to solve them.
The Fisher--KPP equation with nonlocal \cite{lev-trif-shpv1} terms was considered in the
framework of symmetry analysis \cite{olver,bluman}.
Following this approach, we can find
the coefficients that admit a symmetry group of the equation and construct a family of associated particular exact solutions.
Solutions for traveling waves have been found for a reduced Fisher--KPP equation having the only wave independent variable that simplifies
the construction of an exact solution in some special cases \cite{maruvka1,maruvka2,volpert,mei}.

Approximate methods provide additional possibilities in constructing analytical solutions. An example is
the effective particle method \cite{belmonte}, which is in fact the well-known collective coordinate method
(see \cite{bishop1} and references therein). It was originally proposed for Lagrangian equations and then adapted
to the classical (local) Fisher--KPP equation. The method is simple and convenient to use, but, unfortunately, the accuracy of the relevant approximate solutions has never been evaluated or discussed.
To solve the Fisher--KPP equation in a traveling wave study,  homotopy analysis was used (for more details,  see, e.g., \cite{he,shakeel}).

The main objective of the work we present in this paper was to develop a technique which would allow us to apply Maslov complex germ methods
 \cite{maslov1,maslov2,belov1} to the one-dimensional (1D) nonlocal Fisher--KPP equation.

In the general case, for an $n$-dimensional pseudo-differential equation, the construction of semiclassical asymptotic solutions (for both the spectral problem and the Cauchy problem) is based on geometric objects, such as so-called  $k$-dimensional Lagrangian manifolds $\Lambda^k$  with a complex germ $r^n$ ,  where $k<n$. The complex germ $r^n$  is generated by a special set of linearly independent solutions of a linear dynamical system associated with the original equation and called a system in variations.

In this paper, the Maslov complex germ method is adapted to the one-dimensional nonlocal Fisher--KPP equation. To this end, a dynamic Einstein--Ehrenfest system (EE system) is obtained that describes the evolution of the central moments of the solution of the non-local Fisher--KPP equation.
The key point of the proposed approach is that the integrals of the EE system generate the corresponding integrals of the original non-local Fisher--KPP equation.
This allows us to proceed to an auxiliary linear problem associated with the non-local Fisher--KPP equation under consideration. This linear problem is solved by using  symmetry operators, thus giving a way to solve the Cauchy problem for the original non-local Fisher--KPP equation.

The approach developed in this work can also be useful in solving a number of nonlinear problems, including nonlinear cosmological equations (for review see \cite{bamba2012}).

In Section \ref{fkpp}, we consider a one-dimensional nonlocal
Fisher--KPP equation of general form which includes variable
coefficients, local and nonlocal convection, and nonlocal
competitive losses and discuss some basic ideas of the semiclassical
approximation    \cite{maslov1,maslov2,belov1}.
We introduce a class
of trajectory concentrated functions in which the semiclassical
asymptotics are constructed. The asymptotic small parameter here is
the diffusion coefficient; therefore, the asymptotics are obtained
in the approximation of weak diffusion and are valid for a
relatively short time. In Section \ref{eesys}, we deduce a dynamical
system describing the evolution of the moments of an unknown
function of the equation and discuss algebraic conditions for
integration constants. Then we perform a semiclassical reduction of
the original problem to an auxiliary linear problem in Section
\ref{semreduc} and construct semiclassical asymptotics to the
associated linear equation in Section \ref{asslinprob}.
In Section
\ref{leadsemiclsol}, using the results obtained in previous
sections, we find a wide family of approximate particular solutions
in explicit form that have the sense of the leading terms of the
desired semiclassical asymptotics for the nonlocal Fisher--KPP
equation. The family of solutions obtained can be a complete set of
solutions.
In Section \ref{quasistat}, we construct asymptotic
solutions representing a perturbation of some exact homogeneous
solution. These asymptotics  differ  from the semiclassical ones and
are valid for large times.
 In Section \ref{concusion}, we give our
conclusion remarks.

\section{The nonlocal Fisher--KPP equation: a semiclassical approximation}
\label{fkpp}

Consider the following version of the nonlocal Fisher--KPP equation:
\begin{align}
\partial_t u(x,t)&=D\partial_x^2u(x,t)+
  a(x,t)u(x,t)-\partial_x\bigg[V_x (x,t)u(x,t)+\nonumber \\
&+\varkappa u(x,t)\int\limits^\infty_{-\infty} W_x(x,y,t)u(y,t)dy \bigg]-
\varkappa u(x,t)\int\limits^\infty_{-\infty}b(x,y,t) u(y,t)dy.
\label{eq0}
\end{align}
Here $a(x,t)$, $b(x,y,t)$, $V(x,t)$, and $W(x,y,t)$ are given infinitely
smooth functions increasing, as $|x|\,,|y|\to\infty$, no
faster than the polynomial; $\varkappa \,(>0)$ is a real nonlinearity
parameter; $\partial_t =\partial /\partial t$, $\partial_x =\partial /\partial x$,
$\partial_x^2 =\partial^2 /\partial x^2$,  $V_x(x,t)=\partial_x V(x,t)$, and
$W_x(x,y,t)=\partial_x W(x,t)$.

The functions  $V(x,t)$ and $W(x,y,t)$ in equation (\ref{eq0}) are the potentials of the
local and the nonlocal average convective velocity, respectively.

The term $-\varkappa u(x,t)\int\limits^\infty_{-\infty}b(x,y,t) u(y,t)dy$ stands for the nonlocal competition losses
and is characterized by an influence function $b(x,y)$ with the range parameter $\gamma$. To simplify the notation, we
do not introduce $\gamma$ in the function $b(x,y)$ explicitly.

Note that the Fisher--KPP equation with a local convection term was investigated by da Cunha et al.
\cite{cunha2009,cunha2011} and that with a nonlocal convection
term by Ei \cite{ei1987,ei1992}  and Clerc et al. \cite{clerc}.
The one-dimensional  Fisher--KPP equation of the form (\ref{eq0})
was also briefly considered by the authors of Ref. \cite{trif-shpv2009},
where a general scheme of the solution construction was presented.
Here we describe the construction of the semiclassical asymptotics in detail.

The diffusion coefficient $D$ is taken as the small asymptotic parameter.
The supposition of weak diffusion  is the
necessary condition that makes it possible to apply the semiclassical
approach to Cauchy problems with ``narrow'' initial density distributions (i.e.
distributions having a small variance at the time zero). For
the Fisher--KPP equation, this problem is similar to the problem of constructing wave packets in quantum
mechanics \cite{bagr-bel-tr1996}.

Based on the complex WKB method, or the Maslov complex germ
theory  \cite{maslov1,maslov2,belov1}, a method have been developed for constructing asymptotic
solutions in the form of localized wave packets similar to the wave packet solutions of the Schr\"{o}dinger
equation in quantum mechanics \cite{bagr-bel-tr1996}.
 The method of semiclassical asymptotics
was also applied to nonlinear equations; in particular, a Cauchy problem was
solved for a Hartree-type equation \cite{bel-tr-shpv2002}.

To construct asymptotic solutions of equation (\ref{eq0}), we introduce a class of functions
$\mathcal{P}_t^D(X(t,D),S(t,D))$.
The functions of this class singularly
depend on the small parameter
$D$,
a real function
$S(t,D)$ (analog to the action in classical mechanics), and a trajectory $x=X(t,D)$.

The functions $S(t,D)$ and $X(t,D)$  are regular in the
parameter $\sqrt{D}$  in a neighborhood of $D=0$ (i.e. they can
be expressed as a power series in $\sqrt{D}$). The functions of the class
$\mathcal{P}_t^D (X(t,D),S(t,D))$ are concentrated, as $D\to 0$, in a
neighborhood of a point moving in the coordinate space along a
curve given by the equation $x=X(t,$0), and are termed
trajectory-concentrated functions.

Define the common element of the class
$\mathcal{P}_t^D (X(t,D),S(t,D))$ as
\begin{equation}
\label{eq1}
\mathcal{P}_t^D \big(X(t,D),S(t,D) \big)=
  \left\{ \Phi :\Phi(x,t,D)=\varphi\left( \frac{\Delta x}{\sqrt D},t,D\right)
\exp \left[ \frac{1}{D}S(t,D)\right]\right\}.
\end{equation}
Here $\Delta x=x-X(t,D)$. The real function $\varphi (\eta,t,D)$
belongs to the Schwarz space $\mathbb{S}$ in the variable
$\eta\in\mathbb{R}^1$, smoothly depends on $t$, and regularly depends
on $\sqrt{D}$ as $D\to 0$. The real functions $S(t,D)$ and $X(t,D)$,
which characterize the
class, are to be determined. We shall use the contracted notation $\mathcal{P}_t^D $  instead
of $\mathcal{P}_t^D (X(t,D),S(t,D))$ if this does not cause misunderstanding.

By analogy with the work described elsewhere \cite{bagr-bel-tr1996,bel-tr-shpv2002}, it can be shown that
for the functions of the class $\mathcal{P}_t^D $ the following asymptotical estimates hold
\begin{equation}
\label{eq2}
\frac{\|\hat{p}^k \Delta x^lu\|}{\|u\|}\sim O(D^{(k+l)/2}),
\quad
\frac{\|\hat T(X(t,D),t)u\|}{\|u\|}\sim O(D),
\end{equation}
where
\begin{equation}
\label{T-oper}
\hat{p}=D\partial_x, \quad  \hat T(X(t,D),t)=D\partial_t+\dot {X}(t,D)D\partial_x -\dot {S}(t,D).
\end{equation}
In particular,
from (\ref{eq2}) it follows that the operators can be estimated as
$\hat {p}\sim \widehat{O}(\sqrt{D})$ and $\Delta x \sim \widehat{O}(\sqrt{D})$. Here $\widehat O(D^\mu)$ is an
operator $\hat F$ such that
\[\frac {\|\hat F\phi\|}{\|\phi\|}=O(D^\mu),\quad \phi\in\mathcal{P}_t^D .\]
The norm $\|u\|$ is understood in the sense of the space $L_2$.

\section{The Einstein--Ehrenfest system}
\label{eesys}

To develop a method for constructing the semiclassical solutions to equation (\ref{eq0})
in the class of trajectory concentrated functions (\ref{eq1}), we
assume the existence  of the moments
\begin{eqnarray}
\label{eq3}
&&\sigma_u(t,D)=\int\limits^\infty_{-\infty}u(x,t,D)dx,\quad
x_u (t,D)=\frac{1}{\sigma_u(t,D)}
\int\limits^\infty_{-\infty} xu(x,t,D)dx, \\
&&\alpha_u^{(k)}(t,D)=
\frac{1}{\sigma_u(t,D)}\int\limits^\infty_{-\infty}
[x-X(t,D)]^k u(x,t,D)dx .\label{eq4ab}
\end{eqnarray}
Here, $\sigma_u$ and $x_u$  are the zero-order and the first-order moment, respectively, and   $\alpha_u^{(k)}, k\geq 2$,
is the $k$th central moment of a solution  $u(x,t,D) \in \mathcal{P}_t^D $.

From (\ref{eq2}) we obtain the corresponding estimates for  $\sigma_u(t)$, $x_u(t)$, and $\alpha_u^{(k)}$, $k\geq 2$:
\begin{equation}
\label{estim-alpha}
\displaystyle\frac{\sigma_u(t)}{||u||}\sim O(1), \quad   x_u(t) \sim O(1), \quad
\alpha_u^{(k)}(t,D)= O(D^{k/2}),
\end{equation}
and without loss of generality, we can assume that $\sigma_u(t,D)=O(1)$ and $x_u(t,D)=O(1)$.

Let us deduce a dynamical system describing the evolution of the moments (\ref{eq3}),
(\ref{eq4ab}). This is the bicharacteristic system for equation (\ref{eq0}).
The correspondence between the dynamical system and  equation (\ref{eq0}) is similar to that between classical and
quantum mechanics equations. In this context, the system for (\ref{eq3}),
(\ref{eq4ab}) plays the part of the Hamiltonian equations of
classical mechanics and equation (\ref{eq0}) corresponds to the quantum-mechanical Schr\"{o}dinger equation \cite{maslov1,maslov2,bagr-bel-tr1996}.

Let
\begin{gather}
a(x,t)=\sum_{k=0}^{\infty}\frac{1}{k!}a_k(t,D)\Delta x^{k}, \quad
b(x,y,t)=\sum_{k,l=0}^{\infty}\frac{1}{k!l!}b_{k,l}(t,D)\Delta x^{k} \Delta y^{l},\label{expan2}\\
V_x(x,t)=\sum_{k=0}^{\infty}\frac{1}{k!}V_k(t,D)\Delta x^{k}, \quad
W_x(x,y,t)=\sum_{k,l=0}^{\infty}\frac{1}{k!l!}W_{k,l}(t,D)\Delta
x^{k} \Delta y^{l} \label{expan4}
\end{gather}
be formal power series of the functions $a(x,t)$,  $b(x,y,t)$, $V_x(x,t)$, and $W_x(x,y,t)$  in equation (\ref{eq0}) in a neighborhood of $X(t,D)$.
Here, $\Delta y=y-X(t,D)$,
\begin{eqnarray}
&& a_k(t,D)=\frac {\partial^k a(x,t)}{\partial x^k}\bigg|_{x=X(t,D)} , \label{expan1-1}\\
&& b_{k,l}(t,D)=\frac {\partial^{k+l} b(x,y,t)}{\partial x^k\partial y^l}
\bigg|_{\begin{subarray}{l}x=X(t,D)\\
y=X(t,D)\end{subarray}},
\label{expan2-1}\\
&& V_k(t,D)=\frac {\partial^k V_x(x,t)}{\partial x^k}\bigg|_{x=X(t,D)} , \label{expan3-1}\\
&& W_{k,l}(t,D)=\frac {\partial^{k+l} W_x(x,y,t)}{\partial x^k\partial y^l}
\bigg|_{\begin{subarray}{l} x=X(t,D)\\
y=X(t,D)\end{subarray}}.
\label{expan4-1}
\end{eqnarray}
Below we use for  convenience  the shorthand notations $a_k(t,D)=a_k(t)$,  $b_{k,l}(t,D)=b_{k,l}(t)$, $V_k(t,D)=V_k(t)$,
and $W_{k,l}(t,D)=W_{k,l}(t) $.

We impose the technical condition
\begin{equation}
\label{xu}
x_u(t,D)= X(t,D)
\end{equation}
on the functional parameter $X(t,D)$ of the functions $u(x,t,D)$ belonging to the class $\mathcal{P}_t^D$ of trajectory concentrated functions (\ref{eq1}). Then from (\ref{eq3}), (\ref{eq4ab}), and  $\Delta x=x-X(t,D)=x-x_u(t,D)$ it follows that $\alpha_u^{(0)}(t,D)=1$ and $\alpha_u^{(1)}(t,D)=0$.

Differentiating (\ref{eq3}) and (\ref{eq4ab}) with respect to $t$ and substituting $\partial_t u(x,t)$
from (\ref{eq0}) and the expansions (\ref{expan2}) and (\ref{expan4}) in the resulting relations, we obtain the following system of evolution equations for the moments $\sigma_u=\sigma_u(t,D)$, $x_u=x_u(t,D)$, and $\alpha_u^{(k)}=\alpha_u^{(k)}(t,D)$:
\begin{align}
\dot\sigma_u&=\sigma_u \sum\limits^\infty_{k=0}\frac{1}{k!}
 \bigg[ a_k(t) \alpha_u^{(k)}-\sigma_u\varkappa
 \sum_{l=0}^\infty\frac{1}{l!} b_{k,l}(t) \alpha_u^{(k)}\alpha_u^{(l)}\bigg];\label{eq-sigma}\\
\dot x_u&=\sum_{k=0}^\infty \frac{1}{k!}\alpha_u^{(k+1)}\bigg[
  a_k(t)-\sigma_u\varkappa \sum_{l=0}^\infty \frac{1}{l!} b_{k,l}(t)\alpha_u^{(l)}\bigg]+ \nonumber\\
&+\sum_{k=0}^\infty \frac{1}{k!}\alpha_u^{(k)}\bigg[
  V_k(t)+\sigma_u\varkappa \sum_{l=0}^\infty \frac{1}{l!}  W_{k,l}(t)\alpha_u^{(l)}\bigg];\label{eq4}
\\
\dot\alpha_u^{(n)}&=Dn(n-1)\alpha^{(n-2)}_u+\sum_{k=0}^{\infty}\frac{1}{k!}\bigg[
  nV_k(t)\big(\alpha_u^{(k+n-1)}-\alpha_u^{(k)}\alpha_u^{(n-1)}\big)+\nonumber\\
&+a_k(t)\big(\alpha_u^{(k+n)}- \alpha_u^{(k)}\alpha_u^{(n)}-n\alpha_u^{(k+1)}\alpha_u^{(n-1)}\big)\bigg]+\nonumber\\
&+\varkappa
  \sigma_u\sum\limits_{k,l=0}^{\infty}\frac{1}{k!l!}\bigg[
  nW_{k,l}(t)\alpha^{(l)}_u\big(\alpha^{(k+n-1)}_u- \alpha^{(k)}_u\alpha^{(n-1)}_u \big)+ \nonumber\\
&+b_{k,l}(t)\alpha_u^{(l)}\big(- \alpha^{(n+k)}_u+ n \alpha^{(k+1)}_u\alpha^{(n-1)}_u+\alpha^{(k)}_u\alpha^{(n)}_u \big)\bigg].
\label{eq5}
\end{align}
Here $\dot\sigma_u=\displaystyle\frac{d\sigma_u}{dt}$, $\dot x_u=\displaystyle\frac{d x_u}{dt}$, $\dot\alpha_u^{(n)}=\displaystyle\frac{d \alpha_u^{(n)}}{dt}$.

Introducing the aggregate vector of the moments
\begin{equation}
\Theta_u(t, D)=(\sigma_u(t,D), x_u(t,D), \alpha_u^{(2)}(t,D),\alpha_u^{(3)}(t,D) ,  \dots     ),
\label{vector-theta}
\end{equation}
we can rewrite system (\ref{eq-sigma})--(\ref{eq5}) in shorthand form as
\begin{equation}
\dot \Theta_u(t, D)= \Gamma (\Theta_u(t, D), t, D).
\label{system-theta}
\end{equation}
The right-hand side of (\ref{system-theta}),
\begin{equation}
\begin{array}{l}
\Gamma (\Theta_u(t, D),  t, D)=(f(\Theta_u(t, D), t, D), g(\Theta_u(t, D), t, D),\nonumber\\
h^{(2)}(\Theta_u(t, D), t, D), h^{(3)}(\Theta_u(t, D), t, D),  \ldots),\end{array}
\label{system-gamma}
\end{equation}
consists of the functions $f,g, h^{(2)},  h^{(3)}, \dots $  corresponding to the moments $\sigma_u, x_u , \alpha_u^{(2)},  \dots $
the form of which is  obvious from (\ref{eq-sigma}), (\ref{eq4}), and (\ref{eq5}).

The last system (\ref{system-theta}) formally consists of an infinite number of equations. In view of the estimates
(\ref{estim-alpha}), we can reduce this system to a finite number of equations accurate to $O(D^{(M+1)/2})$, conserving
the moments $\alpha_u^{(k)}$ of orders $ k\leq M$ with the order $M$ fixed. The reduced system written in shorthand form like
(\ref{system-theta}) reads
\begin{equation}
\dot \Theta_u^{(M)}(t, D)= \Gamma^{(M)} (\Theta_u^{(M)}(t, D), t, D) + O(D^{(M+1)/2}).
\label{system-theta-r}
\end{equation}
Here we  use the notation
\begin{eqnarray}
&&\hspace*{-1cm}\Theta_u^{(M)}(t, D)=\big(\sigma_u(t,D), x_u(t,D), \alpha_u^{(2)}(t,D),\alpha_u^{(3)}(t,D), \dots ,\alpha_u^{(M)}(t, D) \big),
\label{vector-theta-r}\\
&&\hspace*{-1cm}\Gamma^{(M)} (\Theta_u^{(M)}(t, D), t, D)=\Gamma^{(M)} =\big(f^{(M)}, g^{(M)},
h^{(2,M)}, h^{(3,M)}, \dots,
h^{(M,M)}\big).\label{system-gamma-1}
\end{eqnarray}
In the functions $f^{(M)}, g^{(M)}, h^{(2,M)}, \dots , h^{(M,M)}$  of system (\ref{system-theta-r}),
we conserve the terms of orders not higher than $O(D^{(M/2)})$, using the estimate
$$ \alpha_u^{(k_1)}\alpha_u^{(k_2)} \ldots \alpha_u^{(k_s)} = O(D^{(k_1+k_2+ \dots +k_s)/2})
$$
directly following from (\ref{estim-alpha}).

Consider a system of ordinary differential equations
\begin{equation}
\dot \Theta(t)= \Gamma^{(M)} (\Theta(t); t, D)
\label{eq7a}
\end{equation}
for an aggregate vector
\begin{equation}
\Theta(t)=(\sigma(t), x(t), \alpha^{(2)}(t),\alpha^{(3)}(t),   ...  , \alpha^{(M)}(t)),
\label{vector-theta11}
\end{equation}
where $\Gamma^{(M)} (\Theta(t), t, D)$ is taken from (\ref{system-gamma-1}).

We call equations (\ref{eq7a}) the \textit{Einstein--Ehrenfest}  (EE) \textit{system} of order $M $ for the
Fisher--KPP equation (\ref{eq0}).

The EE system (\ref{eq7a}) is a closed system of evolution equations (i.e. the number of unknown functions is equal to the number of equations),
and hence there exists a unique solution of the Cauchy problem for (\ref{eq7a}).

Suppose that equation (\ref{eq0}) has an exact solution $u(x,t, D)$ (or a solution different from $u(x,t, D)$ by
$O(D^\infty)$) in the class of trajectory-concentrated functions (\ref{eq1})
with the initial condition
\begin{equation}
\label{init-cond}
 u(x,t, D)|_{t=0}=\varphi (x,D)\in\mathcal{P}_0^D.
\end{equation}
Specify the initial condition for system (\ref{eq7a}) as
\begin{equation}
\label{eq7-init}
\Theta(t)|_{t=0}=\mathfrak{g}_\varphi^{(M)}=(\sigma _\varphi ^{(M)}, x_\varphi ^{(M)},
\alpha _\varphi ^{(2,M)},\alpha _\varphi ^{(3,M)} ,\ldots ,\alpha _\varphi^{(M,M)} ),
\end{equation}
where
\begin{equation}
\label{eq6}
\begin{gathered}
\sigma|_{t=0}=\sigma_\varphi^{(M)} =\int\limits_{-\infty}^{\infty}\varphi(x,D)dx,\quad x|_{t=0}=x_\varphi^{(M)}=\frac{1}{\sigma_\varphi}\int\limits_{-\infty}^\infty x\varphi(x,D)dx, \\
\alpha^{(k)}|_{t=0}=\alpha_\varphi^{(k, M)}=\frac{1}{\sigma_\varphi}\int\limits_{-\infty}^\infty (x-x_\varphi)^k
\varphi(x,D)dx, \quad k= 2,  \dots , M.\end{gathered}
\end{equation}
Denote by
\begin{equation}
\label{eq7}
\mathfrak{g}_\varphi^{(M)} (t)=(\sigma _\varphi ^{(M)}(t), x_\varphi ^{(M)} (t),
\alpha _\varphi ^{(2,M)} (t),\ldots ,\alpha _\varphi^{(M,M)} (t))
\end{equation}
the solution of the Cauchy problem for system (\ref{eq7a})  with  the initial condition (\ref{eq7-init}), (\ref{eq6})
and by
\begin{equation}
\label{eq10a}
\mathfrak{g}^{(M)}(t,{\bf C})=(\sigma^{(M)}(t,{\bf C}), x^{(M)}(t,{\bf C}),\alpha^{(2,M)}(t,{\bf C}),
\ldots ,\alpha^{(M,M)}(t,{\bf C}))
\end{equation}
the general solution of system (\ref{eq7a}), where $\bf C$ is a set of integration constants.
We use the shorthand notations $\sigma_\varphi^{(M)}(t)$, $x_\varphi^{(M)} (t)$, and $\alpha_\varphi^{(k,M)}(t)$
for $\sigma_\varphi ^{(M)}(t,D)$, $x_\varphi ^{(M)} (t,D)$, and $\alpha_\varphi^{(k,M)} (t,D)$,  respectively.

The main point in constructing the semiclassical solutions to equation (\ref{eq0}) is that the
solution (\ref{eq10a})  of the the  Einstein--Ehrenfest  system (\ref{eq7a}) can be
obtained independently the solution of the Fisher--KPP equation (\ref{eq0}).

Denote by ${\bf C}_\varphi$ the set of integration constants  ${\bf C}$ in (\ref{eq10a}) determined
by the {\it algebraic} condition
\begin{equation}
\label{algeb-cond-1}
\mathfrak{g}^{(M)}(0,{\bf C})=\mathfrak{g}^{(M)}_\varphi;
\end{equation}
that is, we have $ \mathfrak{g}^{(M)}(0,{\bf C}_\varphi)=\mathfrak{g}^{(M)}_\varphi$. Obviously,  ${\bf C}_\varphi$
is a (vector) functional of $\varphi$.

We do not discuss here the purely algebraic solvability problem  for equation (\ref{algeb-cond-1}) that assumes the existence of ${\bf C}_\varphi$.
In particular cases, ${\bf C}_\varphi$ can be found explicitly (see Section \ref{leadsemiclsol} below).

By virtue of the uniqueness of the solution of the Cauchy problem (\ref{eq7a}), (\ref{eq7-init}) and in view of (\ref{algeb-cond-1}), we have
\begin{equation}
\label{cauchy-EE-1}
\mathfrak{g}^{(M)}(t,{\bf C}_\varphi)=\mathfrak{g}^{(M)}_\varphi (t).
\end{equation}

Consider the aggregate vector
\begin{equation}
\label{cauchy-EE-2}
\mathfrak{g}^{(M)}_u(t)=(\sigma_u(t), x_u(t), \alpha^{(2)}_u(t),  \dots , \alpha^{(M)}_u(t) )
\end{equation}
of the moments $\sigma_u(t), x_u(t), \alpha^{(k)}_u(t)$, $2\leq k\leq M$,  for the solution $u(x,t,D)$ given by (\ref{eq3}) and (\ref{eq4ab}).
Obviously, $\mathfrak{g}^{(M)}_u(t)$ satisfies the EE system (\ref{eq7a}) accurate to $O(D^{(M+1)/2})$:
\begin{equation}
\label{cauchy-EE-4}
\mathfrak{g}^{(M)}_u(t)=\mathfrak{g}^{(M)}_\varphi (t)+ O(D^{(M+1)/2})=\mathfrak{g}^{(M)}(t,{\bf C}_\varphi)+ O(D^{(M+1)/2}),
\end{equation}
and the initial condition
\begin{equation}
\label{cauchy-EE-3}
\mathfrak{g}^{(M)}_u(0)=\mathfrak{g}^{(M)}_\varphi .
\end{equation}
Denote by ${\bf C}_u(t)$ the integration constants ${\bf C}$ in (\ref{eq10a}) obeying the algebraic condition
 (\ref{algeb-cond-1}), i.e.
\begin{equation}
\label{algeb-cond-2}
\mathfrak{g}^{(M)}(t,{\bf C_u(t)})=\mathfrak{g}^{(M)}_u (t).
\end{equation}
Comparing (\ref{algeb-cond-2}) with (\ref{cauchy-EE-4}), we can conclude that
\begin{equation}
\label{algeb-cond-3}
{\bf C}_u(t) = {\bf C}_\varphi + O(D^{(M+1)/2}) .
\end{equation}
This implies that  the  functionals ${\bf C}_u(t)$ provide a set of integrals for the Fisher--KPP equation (\ref{eq0}).

In view of the estimates (\ref{estim-alpha}) and equations (\ref{cauchy-EE-1}),
(\ref{algeb-cond-2}) for $\varphi \in \mathcal{P}_t^D$, we have
\begin{equation}
\label{estim-alpha2}
\frac{\sigma_\varphi (t)}{||\varphi ||}\sim O(1),  \quad   x_\varphi (t) \sim O(1), \quad
\alpha_\varphi^{(k)}(t,D)= O(D^{k/2}).
\end{equation}

To satisfy the estimates (\ref{estim-alpha}), (\ref{estim-alpha2}) in constructing the semiclassical asymptotics
to equation (\ref{eq0}) with the use of the general solution (\ref{eq10a}), we assume that
\begin{eqnarray}
\label{estim-alpha3}
&\sigma^{(M)} (t,\mathbf{C}) \sim O(1), \quad   x^{(M)} (t, {\bf C}) \sim O(1), \\
& \alpha^{(k,M)}(t, {\bf C})= O(D^{k/2}), \quad
\mbox{and} \quad \Delta x=x- x^{(M)}(t, {\bf C})\sim \widehat{O}(\sqrt{D}).\nonumber
\end{eqnarray}

System (\ref{system-theta-r}) comprises information relevant to the solution  $u(x,t,D)$. In particular, the
localization properties of the solution are characterized by the moments  (\ref{eq3}), (\ref{eq4ab}).

As a simple illustrative example, consider system (\ref{eq7a}), restricting the consideration to the terms written accurate to $O(D^{(3/2)})$. To do this, we conserve the terms of orders not higher than  $O(D)$ in (\ref{vector-theta-r}) and write system
(\ref{system-theta-r}) for the  case $M=2$ as
\begin{align}
\dot\sigma_u&=\sigma_u \Bigl[ a(t)+\frac{1}{2}a_2(t)\alpha_u^{(2)} \Bigr]-
\sigma_u^2\varkappa \Bigl[ b(t)+ \frac{1}{2}\big(b_{0,2}(t)+b_{2,0}(t)\big)\alpha_u^{(2)}\Bigr],
\label{sys-111} \\
\dot x_u&=V(t) +\varkappa\sigma_u W(t)+ \Bigl[ a_1(t)+\frac{1}{2}V_2(t)  \nonumber\\
&+\varkappa\sigma_u\Bigl( \frac{1}{2} (W_{0,2}(t)+W_{2,0}(t)) -b_{1,0}(t)\Bigr)\Bigr]\alpha_u^{(2)},
\label{sys-222} \\
\dot\alpha_u^{(2)}&=2D + 2V_1(t) \alpha_u^{(2)}+2\varkappa \sigma_u W_{1,0}(t) \alpha_u^{(2)} .
\label{sys-333}
\end{align}

In particular, for  $a(x,t)=a={\rm const}$, $V(x,t)=W(x,y,t)=0$, $b(x,y,t)=b(x-y)$, and $d b(x)/dx|_{x=0}=0$,
we get equations (\ref{sys-111})--(\ref{sys-333}) in simpler form:
\begin{eqnarray}
&&\dot\sigma_u=a\sigma_u - \varkappa \sigma_u^2 (b+\beta \alpha_u^{(2)}),
\label{sys-112} \\
&&\dot x_u= 0,
\label{sys-223} \\
&& \dot\alpha_u^{(2)}=2D ,
\label{sys-334}
\end{eqnarray}
where $b=b(0)$ and $\beta= d^2 b(x)/dx^2|_{x=0}$.

Equation (\ref{sys-223}) gives $x_u(t,D)=x_u(0,D)$, i.e.  the solution $u(x,t,D)$ in (\ref{eq3}),
(\ref{eq4ab}) is localized in a neighborhood of a fixed point $x_u(0,D)$.

From equation (\ref{sys-334}) we see that $\alpha_u^{(2)}(t,D)=2Dt +\alpha_u^{(2)}(0,D)$, i.e. the
variance of the function $u(x,t,D)$ is the same as that for a linear diffusion with the diffusion coefficient $D$.

Integrating the generalized Verhulst equation (\ref{sys-112}), we obtain
\begin{equation}
\sigma_u (t,D)=a \left[ e^{-at} \left( \frac{a}{\sigma_u(0,D)}-
\varkappa w \right)+ 2D\varkappa \beta t+\varkappa w \right]^{-1},
\label{sys-444}
\end{equation}
where $w=b+\beta \left(\alpha ^{(2)} (0,D)- 2D/a\right)$, and
 (\ref{sys-223}),  (\ref{sys-334}) yield
\begin{eqnarray}
&& x_u (t,D)=x_u (0,D), \nonumber\\
&& \alpha_u^{(2)}(t, D)=2D t+ \alpha_u^{(2)}(0,D).
\label{sys-4444}
\end{eqnarray}
Note that if $b(x,y,t)$ is an asymmetric function in the sense that 
$\partial_x b(x,y,t)|_{x=y=X(t,D)} $ $= v(t,D)\neq 0$,
then $\dot x_u=-\varkappa v\sigma_u\alpha_u^{(2)}$. This implies that the localization region of the solution $u(x,t,D)$
surrounding a point $x_u(t,D)$ moves in space  due to a nonlocal interaction in the system.

\section{Semiclassical reduction  to an auxiliary  linear problem}
\label{semreduc}

Now we can construct semiclassical asymptotic solutions to the nonlocal Fisher--KPP  equation (\ref{eq0}) using the solutions of the EE system (\ref{eq7a}).

To this end, we
expand the functions $b(x,y,t)$ and $W(x,y,t)$ in equation (\ref{eq0}) in powers of $\Delta y=y-x_u (t,D)$
with the remainder of the order $O(D^{(M+1)/2})$
to obtain
\begin{multline}
\bigg\{-\partial_t + D\partial_x^2  + a(x,t)
-\partial_x \bigg[V_x(x,t) +\varkappa \sigma_u(t,D)\displaystyle\sum\limits^M_{l=0} \frac1{l!} W_l(x,x_u, t)\alpha_u^{(l)}(t,D) \bigg]\\
-\varkappa \sigma_u(t,D)\sum\limits^M_{l=0}\frac1{l!} b_l(x,x_u,t) \alpha_u^{(l)}(t,D) \bigg\}u(x,t)=
O(D^{(M+1)/2}),\label{eq10}
\end{multline}
where
\begin{eqnarray}
&W_l(x,x_u,t)=\dfrac{\partial^l W_x(x,y,t)}{\partial y^l}\Big|_{y=x_u(t,D)},
\label{lin-eq-1} \\
& b_l(x,x_u,t)=\dfrac{\partial^l b(x,y,t)}{\partial y^l}\Big|_{y=x_u(t,D)}.
\label{lin-eq-2}
\end{eqnarray}
The moments $\sigma_u(t,D)$ and $\alpha_u^{(l)}(t,D)$  are given by (\ref{eq3}) and (\ref{eq4ab}), $l=0,1, \dots ,M$.

Substituting  the general  solution (\ref{eq10a})  of the Einstein--Ehrenfest system (\ref{eq7a}) of order
$M$ instead of the corresponding moments $\sigma_u$, $x_u$, $\alpha_u^{(k)}$, $k= 2,\dots,M$, in (\ref{eq10}),
we obtain  a {\it linear} partial differential equation
\begin{eqnarray}
&&\hat L (x,t,  {\bf C})v(x,t,{\bf C})= \nonumber\\
&&\hspace*{-1cm}=\big[-\partial_t +D\partial_x^2 -\partial_x \Lambda^{(M)}(x,t,{\bf C})+A^{(M)}(x,t,{\bf C})
\big]v(x,t,{\bf
C})=O(D^{(M+1)/2}),
\label{eq11}
\end{eqnarray}
where
\begin{equation}
\label{eq12}
\Lambda^{(M)}(x,t,{\bf C})=V_x (x,t)+\varkappa
\sigma^{(M)}(t,{\bf C})\sum^M_{l=0}\frac1{l!} W_l(x,x^{(M)}(t,{\bf C}),t)\alpha^{(l,M)}(t,{\bf C}),
\end{equation}
\begin{equation}
\label{eq13}
A^{(M)}(x,t,{\bf C})=a(x,t)-\varkappa
\sigma^{(M)}(t,{\bf C})\sum^M_{l=0}\frac1{l!}
b_l(x,x^{(M)}(t,{\bf C}),t)
\alpha^{(l,M)}(t,{\bf C}).
\end{equation}
Here, $\alpha^{(0,M)}(t,{\bf C})=1$, $\alpha^{(1,M)}(t,{\bf C})=0$, and the form of $W_l(x,x^{(M)}(t,{\bf C}),t )$ and  $b_l(x,x^{(M)}(t,{\bf C}),t )$ is clear from
(\ref{lin-eq-1}) and (\ref{lin-eq-2}).

Equation (\ref{eq11}) can be considered as a family of linear PDEs parametrized by  the arbitrary integration constants
${\bf C}$.

We shall term equation (\ref{eq11}) with the coefficients (\ref{eq12}), (\ref{eq13}) \textit{the associated linear Fisher--KPP
equation} of order $M$.

Note that, according to (\ref{eq13}) and (\ref{eq-sigma}), we can write
\begin{eqnarray}
&&A^{(M)}(x,t,{\bf C}) - \frac{\dot\sigma^{(M)}(t,{\bf C})}{\sigma^{(M)}(t,{\bf C})}=
\sum\limits^M_{k=0}\frac{1}{k!}
 \bigg[ a_k(t) [(\Delta x)^k-\alpha^{(k,M)}(t,{\bf C})]\nonumber\\
&&\quad -\sigma^{(M)}(t,{\bf C})\varkappa
 \sum_{l=0}^M\frac{1}{l!} b_{k,l}(t) [(\Delta x)^k-\alpha^{(k,M)}(t,{\bf C})]\alpha^{(l,M)}(t,{\bf C})
 \bigg];\label{eq-sigma-1}
\end{eqnarray}
therefore, we have
\begin{eqnarray}
&& A^{(M)}(t,{\bf C})=a(t,{\bf C})-\varkappa
\sigma^{(M)}(t,{\bf C})\sum^M_{l=0}\frac1{l!}
b_l(t,{\bf C})\alpha^{(l,M)}(t,{\bf C})\nonumber\\
&&\quad =\displaystyle\frac{\dot\sigma^{(M)}(t,{\bf C})}{\sigma^{(M)}(t,{\bf C})}- \sum\limits^M_{k=2}\displaystyle\frac{1}{k!}
 \bigg[ a_k(t) \alpha^{(k,M)}(t,{\bf C})\nonumber\\
&&\quad -\sigma^{(M)}(t,{\bf C})\varkappa
 \sum_{l=2}^M\displaystyle\frac{1}{l!} b_{k,l}(t) \alpha^{(k,M)}(t,{\bf C})\alpha^{(l,M)}(t,{\bf C})
 \bigg]+O(D^{(M+1)/2}).\label{eq-sigma-2}
\end{eqnarray}
Here and below we set $A^{(M)}(t,{\bf C})=A^{(M)}(x^{(M)}(t,{\bf C}),t,{\bf C})$; in particular,
\begin{equation}
A^{(M)}(t,{\bf C})=\frac{\dot\sigma^{(M)}(t,{\bf C})}{\sigma^{(M)}(t,{\bf C})}+O(D).\label{eq-sigma-3}
\end{equation}

Let us  reduce the search for a solution to the Cauchy problem for the Fisher--KPP equation (\ref{eq10}) in the class $\mathcal{P}_t^D$
of trajectory concentrated functions accurate to  $O(D^{(M+1)/2})$ to the search for a solution to the corresponding Cauchy problem for the
associated linear equation (\ref{eq11}), using the general solution (\ref{eq10a}) of the EE system (\ref{eq7a}).

To this end, we  consider the Cauchy problems for equations  (\ref{eq10}) and (\ref{eq11}) with the same
initial condition
\begin{equation}
\label{eq14}
u\vert _{t=0}=v\vert_{t=0}=\varphi(x, D),\quad
\varphi(x, D)\in \mathcal{P}_0^D .
\end{equation}

Replace the arbitrary constants $ {\bf C}$ in equation (\ref{eq11}) by  the constants ${\bf C}_\varphi $ found from (\ref{algeb-cond-1}).

In view of  (\ref{cauchy-EE-1}) and (\ref{cauchy-EE-4}), the associated linear equation (\ref{eq11})
transforms into equation (\ref{eq10}) accurate to $O(D^{(M+1)/2})$.

Then, in virtue of the uniqueness of the Cauchy problem, we obtain the key theorem reducing  the construction of semiclassical asymptotics
for the nonlinear equation (\ref{eq0}) to the search for a solution to the linear equation (\ref{eq11}).

\textbf{Theorem 1.}
{\it The solutions of the Cauchy problem  for the nonlinear equation
\eqref{eq10} and of the Cauchy problem for the associated linear
equation \eqref{eq11} with the same initial condition \eqref{eq14}
are related as
\begin{equation}
\label{eq15}
u(x,t)=v^{(M)}(x,t,{\bf C}_\varphi)+O(D^{(M+1)/2}),
\end{equation}
where the constants ${\bf C}_\varphi$ are determined by {\rm(\ref{algeb-cond-1})}}.

Note that the change from the Cauchy problem (\ref{eq10}), (\ref{eq14}) for $ u(x,t)$ to the Cauchy problem (\ref{eq11}),
(\ref{eq14}) for  $v^{(M)}(x,t,{\bf C})$ is not the usual linearization of equation (\ref{eq11}) by expanding its terms in power series in a small parameter.
Equation (\ref{eq11}) has the sense of an auxiliary problem for (\ref{eq10}).

In the next section,  we  construct a solution to the Cauchy problem for the associated linear equation (\ref{eq11})
in terms of semiclassical asymptotics.

\section{Semiclassical asymptotics of the associated linear equation}
\label{asslinprob}

Main idea of the search for a semiclassical solution to the Cauchy problem for the associated linear equation (\ref{eq11}) can be formulated as follows:

Take a function $\varphi(x,D) $ which belongs to the  class $\mathcal{P}_0^D$ of the form (\ref{eq1}) at $t=0$.

Let $v_{\rm ext}(x,t, {\bf C}, D)$ be an exact solution of the Cauchy problem  to equation (\ref{eq11})
with the initial condition $\varphi(x, D)$, and a function $v_{\rm as}(x,t, {\bf C}, D)$
of the class $\mathcal{P}_t^D$ be such that $v_{\rm as}(x,0, {\bf C}, D)=v_{\rm as}(x,0, D)=$ $\varphi(x, D)$.
Denote by  $h^{(\mu)}(x,t, {\bf C}, D)$  the residual of the function
$v_{\rm as}(x,0,  D)$ for the operator $\hat L (x,t, {\bf C})$ in equation  (\ref{eq11}):
\begin{equation}
\label{ale-1}
\hat L (x,t,{\bf C}) v_{\rm as}(x,t, {\bf C}, D)=h^{(\mu)}(x,t, {\bf C}, D).
\end{equation}

According to \cite{maslov1}, the following statement holds:
if $||h^{(\mu)}||=O(D^{\mu})$, then
\begin{equation}
\label{residual}
||v_{\rm ext}-v_{\rm as}||=O(D^{\mu}).
\end{equation}

The estimate  (\ref{residual}) is valid for the functions $v_{\rm ext}$ and  $v_{\rm as}$ belonging the class $\mathcal{P}_t^D$.
As the solution  of the Cauchy problem (\ref{eq11}), (\ref{eq14}) is continuous in $t$, the function $v_{\rm ext}$ belongs the class $\mathcal{P}_t^D$
at least in a certain time interval $[0,T]$, $t\in[0,T]$. Analysis  of the asymptotic solutions for the Fisher--KPP equation (\ref{eq0}) in particular cases (see, e.g., \cite{12c}) shows that the residual $h^{(\mu)}$ can either increase or decrease as  $t\to\infty$.
Therefore, the problem of finding the time interval $[0,T]$ where the semiclassical asymptotics hold requires a separate investigation that goes beyond the scope of this work.

Relation (\ref{residual}) implies that  asymptotic solutions  of the Cauchy problem should be constructed under the condition $\mu >0$.
Accordingly, it suffices to find the leading order term of the asymptotic solutions  in the class $\mathcal{P}_t^D$ accurate to $O(D^{1/2})$.

Now we  construct the semiclassical asymptotic solution of the Cauchy problem for  the associated linear equation (\ref{eq11}).

Expand the coefficients of the operator $\hat L(x,t,{\bf C})$ in terms of powers of $\Delta x=x-x^{(M)} (t,{\bf C})$, where $x^{(M)}(t,{\bf C})$
is a component of the general solution (\ref{eq10a}).
Then  the operator $D\hat L(x,t,{\bf C})$ reads
\begin{eqnarray}
&&D\hat L(x,t,{\bf C})= -\hat T(x^{(M)}(t, {\bf C}),t)+\dot x^{(M)}(t, {\bf C})\hat p-\dot S(t,D)+\hat p^2 \nonumber\\
&&\quad-\hat p\sum_{l=0}^{M}\displaystyle\frac{1}{l!}\Lambda_{(l)}^{(M)}(t,{\bf C})\Delta x^l
+D\sum_{l=0}^{M}\displaystyle\frac{1}{l!} A_{(l)}^{(M)}(t,{\bf C})\Delta x^l+\hat O(D^{l/2 +1}).
\label{L-expan-1}
\end{eqnarray}
Here, $\hat T(x^{(M)}(t, {\bf C}),t)=D\partial_t + \dot x^{(M)}(t, {\bf C})\hat p-\dot S(t,D)$ and $\hat p$ are given by
(\ref{T-oper}),
\begin{equation}
\begin{aligned}
&\Lambda_{(l)}^{(M)}(t,{\bf C})=\dfrac{\partial^l \Lambda ^{(M)}(x,t, {\bf C})}{\partial x^l}|_{x=x^{(M)}(t, {\bf C})},\\
&A_{(l)}^{(M)}(t,{\bf C})=\dfrac{\partial^l A^{(M)}(x,t, {\bf C})}{\partial x^l}|_{x=x^{(M)}(t, {\bf C})}.\end{aligned}
\label{L-expan-2}
\end{equation}

The function $S(t,D)$ is a free parameter of the class (\ref{eq1}).
Choosing $S(t,D)$ such that the terms proportional to the identity operator vanish in the operator (\ref{L-expan-1}), we have
\begin{equation}
\dot S(t,D)=DA^{(M)}(x^{(M)}(t, {\bf C}),t, {\bf C} ).
\label{s-function-1}
\end{equation}

From (\ref{eq-sigma-3}) and  (\ref{s-function-1}) it follows that
\begin{equation}
\exp\Big[\frac1D S(t,D)\Big]=\exp\Big[\int_0^t\frac{\dot\sigma^{(M)}(t,{\bf C})}{\sigma^{(M)}(t,{\bf C})}dt+O(D)\Big]= 
\frac{\sigma^{(M)}(t,{\bf C})}{\sigma^{(M)}(0,{\bf C})}+O(D).\label{eq-sigma-4}
\end{equation}

The evolution equation of the EE system (\ref{eq7a}) for $x^{(M)}(t, {\bf C})$ reads
\begin{align}
\dot x^{(M)}&=\sum_{k=0}^{M}\dfrac{1}{k!}\alpha^{(k+1,M)}\big[ a_k-\sigma^{(M)}\varkappa
\sum_{l=0}^{M}\dfrac{1}{l!}b_{k,l}\alpha^{(l,M)} \big]  \nonumber\\
&+\sum_{k=0}^{M}\dfrac{1}{k!}\alpha^{(k,M)}\big[V_k+\sigma^{(M)}\varkappa \sum_{l=0}^{M}\dfrac{1}{l!}W_{k,l}\alpha^{(l,M)}\big],
\label{L-expan-4}
\end{align}
where $a_k,  b_{k,l},V_k,W_{k,l}$ are given by (\ref{expan1-1})--(\ref{expan4-1}),
$\sigma^{(M)}=\sigma^{(M)}(t,{\bf C})$, $x^{(M)}=x^{(M)}(t,{\bf C})$, and $\alpha^{(k,M)}=\alpha^{(k,M)}(t,{\bf C})$.
Insert $\dot x^{(M)}(t,{\bf C})$ from (\ref{L-expan-4}) in (\ref{L-expan-1}). Using estimates (\ref{eq2}), (\ref{T-oper}), and (\ref{estim-alpha3}),
we present $D\hat L(x,t,{\bf C})$ in (\ref{L-expan-1}) as an expansion:
\begin{equation}
D\hat L(x,t,{\bf C}) = \hat L_{(0)}^{(M)}(t,{\bf C})+ \hat L_{(1)}^{(M)}(t,{\bf C})+\hat L_{(2)}^{(M)}(t,{\bf C}) \dots ,
\label{L-expan-5}
\end{equation}
where  $\hat L_{(0)}^{(M)}(t,{\bf C})\sim \hat O(D) $, $\hat L_{(1)}^{(M)}(t,{\bf C})\sim \hat O(D^{3/2})$, and
$\hat L_{(2)}^{(M)}(t,{\bf C})\sim \hat O(D^{2})$, ... .
After some calculations, we find  $\hat L_{(k)}^{(M)}(t,{\bf C})$, $k=0, 1, 2, \dots $, as
\begin{align}
&\hat L_{(0)}^{(M)}(t,{\bf C})= -\hat T(x^{(M)}(t,{\bf C}),t)+\hat
p^2 -\hat p\Lambda_x^{(M)}(x^{(M)}(t,{\bf C}),t,{\bf C})\Delta x,
\label{L-expan-6}\\
&\hat L_{(1)}^{(M)}(t,{\bf C})=\alpha^{(2,M)} (t, {\bf C})\biggl[
a_x(x^{(M)}(t,{\bf C}), t) \nonumber\\
&\qquad- \sigma^{(M)}(t,{\bf C})\varkappa \sum_{l=0}^{M}\frac{1}{l!}
b_{1,l}\alpha^{(l,M)} (t, {\bf C})+
\frac{1}{2}V_2(t)\nonumber\\
&\qquad+\frac{1}{2}\sigma^{(M)}(t,{\bf C})\varkappa \sum_{l=0}^{M}
\frac{1}{l!}W_{2,l}(t)\alpha^{(l,M)}(t,{\bf C})\biggr]\hat p \nonumber\\
&\qquad -\frac{1}{2}\Lambda^{(M)}_{xx}(x^{(M)}(t,{\bf C}), t, {\bf
C})\hat p\Delta x^2+
DA^{(M)}_x(x^{(M)}(t, {\bf C}), t,{\bf C})\Delta x,\label{L-expan-7}\\
&\hat L_{(2)}^{(M)}(t,{\bf C})=\frac{1}{2}\alpha^{(3,M)} (t, {\bf
C})\bigg[
a_2(x^{(M)}(t,{\bf C}), t)\nonumber\\
&\qquad -\sigma^{(M)}(t,{\bf C})\varkappa \sum_{l=0}^{M}\frac{1}{l!} b_{2,l}\alpha^{(l,M)} (t, {\bf C})\nonumber\\
&\qquad+\frac{1}{3}\biggl(V_3(t)+\sigma^{(M)}(t,{\bf C})\varkappa \sum_{l=0}^{M}\frac{1}{l!}W_{3,l}(t)\alpha^{(l,M)}(t, {\bf C})\biggr)\bigg]\hat p \nonumber\\
&\qquad - \frac{1}{3!}\Lambda^{(M)}_{(3)}(x^{(M)}(t,{\bf C}), t,
{\bf C})\hat p\Delta x^3+
\frac{1}{2}DA^{(M)}_{(2)}(x^{(M)}(t, {\bf C}), t,{\bf C})\Delta x^2,\label{L-expan-8}\\
...&  .................................\nonumber
\end{align}

Note that the operators $\hat L_{(k)}^{(M)}(t,{\bf C})$, $k=0,1,2, \dots $, in (\ref{L-expan-5}) are not uniquely determined.

We call the expansion
$$
D\hat L(x,t,{\bf C}) = \hat{\tilde L}_{(0)}^{(M)}(t,{\bf C})+ \hat{\tilde L}_{(1)}^{(M)}(t,{\bf C})+\hat{\tilde L}_{(2)}^{(M)}(t,{\bf C})+ \dots
$$
to be equivalent to (\ref{L-expan-6}) if
$$
\hat{\tilde L}_{(k)}^{(M)}-\hat{L}_{(k)}^{(M)}\sim \hat O(D^{(k+1)/2}).
$$

We shall seek the trajectory-concentrated solutions $v(x,t,D)\in
\mathcal{P}_t^D$ to the linear associated equation (\ref{eq11}) in the form
\begin{equation}
\label{eq19}
v(x,t,{\bf C},D)=v^{(0)}(x,t,{\bf C},D)+\sqrt{D}v^{(1)}(x,t,
{\bf C},D)+Dv^{(2)}(x,t,{\bf C},D)+\ldots ,
\end{equation}
where $v^{(k)}(x,t,{\bf C},D) $ $\in\mathcal{J}_t^D (x^{(M)}(t,{\bf
C}),S(t,D))$.
The class of functions\\
$\mathcal{J}_t^D (x^{(M)}(t,{\bf
C}),S(t,D))\subset $ $\mathcal{P}_t^D (x^{(M)}(t,{\bf C}),S(t,D))$ is
obtained from $\mathcal{P}_t^D$ if a function $\varphi(\eta,t)$
entering into the definition of the class does not depend on the
parameter $D$.

Substituting the expansion (\ref{eq19}) in the linear  equation (\ref{eq11}), we have
\begin{equation}
D\hat L(x,t,{\bf C})v(x,t,{\bf C},D)=O(D^{(M+1)/2}),
\label{eq333}
\end{equation}
where $D\hat L(x,t,{\bf C})$ is taken in the form  (\ref{L-expan-5}).
Equating the terms having the same estimate in the parameter $D$, in view of the estimates
(\ref{eq2}), (\ref{estim-alpha}), and (\ref{estim-alpha3}),
we obtain the following system of recurrent equations:
\begin{eqnarray}
&& \hat {L}_{(0)}^{(M)}(t,{\bf C})v^{(0)}=0,\label{eq23a}\\
&& \hat {L}_{(0)}^{(M)}(t,{\bf C})v^{(1)}+
\hat{L}_{(1)}^{(M)}(t,{\bf C})v^{(0)}=0, \label{eq23b}\\
&& \hat{L}_{(0)}^{(M)}(t,{\bf C})v^{(2)}+
\hat{L}_{(1)}^{(M)} (t,{\bf C})v^{(1)}+
\hat {L}_{(2)}^{(M)}(t,{\bf C})v^{(0)}=0, \label{eq23c}\\
&& \ldots \ldots \ldots \ldots \ldots \ldots \ldots \ldots \ldots \ldots
\ldots\nonumber
\end{eqnarray}

Thus, solving system (\ref{eq23a})--(\ref{eq23c}) to a given accuracy $O(D^{(M +1)/2})$, we can obtain
the desired asymptotic solution of the nonlocal  Fisher--KPP equation (\ref{eq0}) corresponding to  (\ref{eq15}).

\section{The leading order terms  of semiclassical solutions }
\label{leadsemiclsol}

We construct here an explicit family of particular solutions $v^{(0)}(x,t,{\bf C})$ for  equation
(\ref{eq23a}). Such solutions are widely used in the Maslov complex germ theory  \cite{maslov1,maslov2,belov1}
 and form a convenient basis for
solutions of the associated linear equation (\ref{eq23a}).
With the help of the solutions obtained, we find the corresponding  leading order terms of the semiclassical asymptotic expansions for the nonlocal Fisher--KPP equation.

Let us seek a particular solution to equation (\ref{eq23a}) in the form
\begin{equation}
\label{eq20}
v_0^{(0)}(x,t,{\bf C})=N_D\exp \left\{ \frac{1}{D}\left(
S(t,D)+\frac{1}{2}Q(t)\Delta x^2\right) + \varphi (t,D)\right\},
\end{equation}
where $N_D$ is a constant, the function $Q(t)$ is to be determined, and
$S(t,D)$ is determined by (\ref{s-function-1}).
Substituting (\ref{eq20}) in (\ref{eq23a}) and equating the coefficients of
the terms with identical powers of $\Delta x$, we obtain the following system of equations:
\begin{equation}
\label{eq21}
\dot \varphi + \Lambda _x^{(M)} (x^{(M)}(t,{\bf C}),t,{\bf C})-Q=0,
\end{equation}
\begin{equation}
\label{eq22}
\frac{1}{2}\dot{Q}- Q^2 + \Lambda_x^{(M)}(x^{(M)}(t,{\bf C}),t,{\bf C})Q=0.
\end{equation}
From (\ref{eq21}) it follows that
\begin{equation}
\label{eq23}
\varphi (t,D)=\int\limits^t_0[-\Lambda_x^{(M)}(x^{(M)}(t,{\bf C}),t,{\bf C})+Q(t)]dt.
\end{equation}

The linear equation (\ref{eq23a}) is quadratic in $x$ and $\partial_x$, and therefore (see, e.g., \cite{maslov1}) it
admits a symmetry operator linear in $x$ and $\partial_x$. Let us seek this operator in the form
\begin{equation}
\hat{a}(t,{\bf C})=N_a [Z(t,{\bf C})\hat{p}-W(t,{\bf C})\Delta x],
\label{28a}
\end{equation}
where $N_a$  is a constant factor. The functions $Z(t,{\bf C})$ and
$W(t,{\bf C})$ are determined from the condition that the
commutator should be equal to zero:
\begin{equation}
\label{eq24}
[\hat {L}_0 (t,{\bf C}),\hat{a}(t,{\bf C})]=0.
\end{equation}
Here $\hat {L}_0=-\hat T+\hat p^2- \Lambda_x^{(M)}\hat p \Delta x $,
where $\hat T=D\partial_t+\dot x^{(M)}\hat p-\dot S$ and $\hat p=D\partial_x$.

Relation (\ref{eq24}) will be hold if the functions
$Z(t,{\bf C})$ and $W(t,{\bf C})$ satisfy the  equations
\begin{equation}
\label{eq25}
\begin{cases}
\dot{W}=-\Lambda_x^{(M)} (x^{(M)}(t,{\bf C}),t)W,\\
\dot {Z}=-2W+\Lambda_x^{(M)} (x^{(M)}(t,{\bf C}),t)Z. \end{cases}
\end{equation}
In addition,
\begin{equation}
\label{eq26}
Q(t,{\bf C})=\frac{W(t,{\bf C})}{Z(t,{\bf C})}
\end{equation}
is the solution of equation (\ref{eq22}), as one can check directly.
Equations (\ref{eq25}) constitute  the system in variations  \cite{maslov2} to
equation (\ref{eq22}).

From (\ref{eq25}) and (\ref{eq26}) we obtain
\begin{equation}
\label{eq25-1}
\begin{cases}
\Lambda_x^{(M)} (x^{(M)}(t,{\bf C}),t)=-\dfrac{\dot{W}}{W},\\
Q-\dfrac12\Lambda_x^{(M)} (x^{(M)}(t,{\bf C}),t)=-\dfrac12\dfrac{\dot {Z}}{Z}.\end{cases}
\end{equation}

Then (\ref{eq23}) can be written as
\begin{align}
\exp\big[\varphi (t,D)\big]&=\exp\biggl[\int\limits^t_0[-\Lambda_x^{(M)}(x^{(M)}(t,{\bf C}),t,{\bf C})+Q(t)]dt\biggr]=\nonumber\\
&=\exp\biggl[\frac12\int_0^t\Big(\frac{\dot{W}(t,{\bf C})}{W(t,{\bf C})}-\frac{\dot {Z}(t,{\bf C})}{Z(t,{\bf C})}\Big)dt\biggr]=\sqrt{\bigg|\frac{Z(0,{\bf C})W(t,{\bf C})}
{Z(t,{\bf C})W(0,{\bf C})}\bigg|}\label{eq-sigma-5}
\end{align}
and, in view of (\ref{eq-sigma-4}), solution (\ref{eq20}) can be presented as
\begin{equation}
\label{eq27} v_0^{(0)} (x,t,{\bf C})=N_D \sqrt{\bigg|\frac{Z(0,{\bf C})W(t,{\bf
C})} {Z(t,{\bf C})W(0,{\bf C})}\bigg|} \frac{\sigma^{(M)}(t,{\bf C})}
{\sigma^{(M)}(0,{\bf C})}\exp\Big[\frac{1}{2D}\frac{W(t,{\bf C})}
{Z(t,{\bf C})}\Delta x^2\Big].
\end{equation}
Using the solutions of the system in variations (\ref{eq25}),
we compose a two-dimensional vector $a(t)=a(t,{\bf C})=
(W(t,{\bf C}),Z(t,{\bf C}))$ and designate by $a^{(+)}(t,{\bf
C})=(W^{(+)}(t),Z^{(+)}(t))$ and $a^{(-)}(t,{\bf C})=
(W^{(-)}(t),Z^{(-)}(t))$ the solutions of the system in variations
that satisfy the conditions
\begin{equation}
W^{(+)}(0)=b,\quad W^{(-)}(0)=-b,\quad b>0,\quad Z^{(\pm )}(0)=1.
\end{equation}
Note the fact essential to the subsequent considerations: For the
solutions of the system in variations, the
following skew-scalar product is preserved:
\begin{equation}
\label{eq28}
\{a^{(-)}(t), a^{(+)}\}=W^{(+)}(t)Z^{(-)}(t)-Z^{(+)}(t)W^{(-)}(t)=2b.
\end{equation}
The pair $\big(\Lambda^0_t, r^1_t \big)$  is a manifold  $\Lambda^0_t$  with a germ $r^1_t$. In our case, the germ $r^1_t$  is chosen to be real as we seek real solutions to the non-local Fisher--KPP equation (\ref{eq0}).
Here, the  0-dimensional manifold  $\Lambda^0_t$  is a point in the phase space of equation (\ref{sys-222}) which is given at every time $t$ by the relation $x=X(t)$. The germ $r^1_t$  is a linear 1-dimensional space with a basis vector $a^{(-)}(t)$.

The vectors $a^{(\pm)}(t)=(Z^\pm(t),W^{(\pm)}(t))$ are associated with
two symmetry operators of the form (\ref{28a}): $\hat {a}^{(-)}(t,{\bf C})$
and $\hat{a}^{(+)}(t,{\bf C})$, which, following the terminology of
quantum mechanics, can be naturally termed the operators of ``annihilation''
and ``creation'' of states, respectively. Using formula (\ref{eq27}), we can
construct two solutions to equation (\ref{eq23a}): $v_0^{-(0)}(x,t,{\bf C})$
and $v_0^{+(0)}(x,t,{\bf C})$. In what follows, any of the functions
$v_0^{-(0)}(x,t,{\bf C})$ and $v_0^{+(0)}(x,t,{\bf C})$ will
be required to us. For definiteness, we choose the function
$v_0^{-(0)}(x,t,{\bf C})$, and hereinafter we shall omit the superscript
``$-$'' where it does not cause misunderstanding:
\begin{equation}
\label{eq29}
v_0^{(0)} (x,t,{\bf C})=N_D\sqrt{-\frac{W^{(-)}(t)}{bZ^{(-)}(t)}}
\frac{\sigma^{(M)}(t,{\bf C})}{\sigma^{(M)}(0,{\bf C})}\exp\Big[
\frac{1}{2D}\frac{W^{(-)}(t)}{Z^{(-)}(t)}\Delta x^2\Big].
\end{equation}
We also assume that $W^{(-)}(t)Z^{(-)}(t)<0$.

Normalizing the operators $\hat{a}^{(\pm)}(t,{\bf C})$ by the
condition
$$
[\hat {a}^{(-)},\hat {a}^{(+)}]=1,
$$
we then have
\begin{equation}
\label{eq30}
\hat {a}^{(+)}(t,{\bf C})=\frac{1}{\sqrt{2bD}}
[Z^{(+)}(t)D\partial_x -W^{(+)}(t)\Delta x],
\end{equation}
\begin{equation}
\label{eq30aa}
\hat {a}^{(-)}(t,{\bf C})=-\frac{1}{\sqrt{2bD}}
[Z^{(-)}(t)D\partial_x -W^{(-)}(t)\Delta x].
\end{equation}
Note that the function given by (\ref{eq29}) describes a ``vacuum''
state for the operator $\hat {a}^{(-)}(t,{\bf C})$:
\begin{equation}
\label{eq31}
\hat{a}^{(-)}(t,{\bf C})v_0^{(0)}(x,t,{\bf C})=0.
\end{equation}
Set the normalizing factor $N_D $ by
$$
N_D=\sqrt[4]{\frac{b}{\pi D}}.
$$
Let us specify a set of functions resulting from sequential actions of the
creation operators on the vacuum state:
\begin{equation}
\label{eq32}
v_n^{(0)}(x,t,{\bf C})=\frac{1}{\sqrt{n!}}
(\hat{a}^{(+)}(t,{\bf C}))^nv_0^{(0)}(x,t,{\bf C}).
\end{equation}
Using the Rodrigues formula for Hermite polynomials \cite{Beitman2}
$$
H_n(\xi)=(-1)^ne^{\xi^2}(e^{-\xi^2})^{(n)},
$$
we write
\begin{equation}
\label{eq33}
v_n^{(0)}(x,t,{\bf C})=\frac{(-1)^n}{\sqrt{2^nn!}}
\Big[\frac{Z^{(+)}(t)}{Z^{(-)}(t)}\Big]^{n/2}H_n
\left(\sqrt{\frac{b}{DZ^{(-)}(t)Z^{(+)}(t)}}\Delta x\right)
v_0^{(0)}(x,t,{\bf C}).
\end{equation}

From the definitions of the operator $\hat{a}^{(+)}(t,{\bf C})$ and
function $v_0^{(0)}(x,t,{\bf C})$ it follows that the functions given
by (\ref{eq32}) will
be a solution of equation (\ref{eq23a}) for an arbitrary natural $n$.

Let us consider an auxiliary set of functions:
\begin{equation}
\label{eq34}
w_n^{(0)}(x,t,{\bf C})=\frac{(-1)^n}{\sqrt{2^nn!}}
\Big[\frac{Z^{(-)}(t)}{Z^{(+)}(t)}\Big]^{n/2}H_n
\left(\sqrt{\frac{b}{DZ^{(+)}(t)Z^{(-)}(t)}}\Delta x\right)
w_0^{(0)}(x,t,{\bf C}),
\end{equation}
where
\begin{equation}
\label{eq35}
w_0^{(0)}(x,t,{\bf C})=
\frac{N_D\sigma(0)}{\sigma^{(M)}(t,{\bf C})}\sqrt{\frac{-b}{W^{(-)}(t)Z^{(+)}(t)}}\exp\Big[
-\frac{1}{2D}\frac{W^{(+)}(t)}{Z^{(+)}(t)}\Delta x^2\Big].
\end{equation}
By analogy with the terminology of quantum mechanics, we shall term the functions given
by (\ref{eq32}) \textit{the semiclassical trajectory-coherent solutions}
of the Fisher--KPP equation (\ref{eq0}).

The  functions $\{v_n^{(0)}(x,t,{\bf C})\}_{n=0}^\infty$
of the form (\ref{eq33}) and the functions  
$\{w_m^{(0)}(x,t,{\bf C})\}_{m=0}^\infty$
of the form (\ref{eq34}) are biorthogonal:
\begin{equation}
\label{eq37}
\int\limits^{+\infty}_{-\infty}v_n^{(0)}(x,t,{\bf C} )w_m^{(0)}
(x,t,{\bf C})dx=\delta_{nm},
\end{equation}
and they satisfy the completeness condition
\begin{equation}
\sum^\infty_{n=0}v_n^{(0)}(x,t,{\bf C})w_n^{(0)}
(y,t,{\bf C})=\delta(x-y).
\label{44a}
\end{equation}

According to Theorem 1  (see (\ref{eq15})), the asymptotic
solutions of the nonlinear equation (\ref{eq0}), accurate to  $O(D^{3/2})$, are
\begin{equation}
\label{solution-family}
u_n (x,t)=v_n^{(0)}(x,t,{\bf C}_n).
\end{equation}
Here, the constants  ${\bf C}_n$  are determined by the equation
$\mathfrak{g}^{(2)}(t,{\bf C}_n)=\mathfrak{g}_{\varphi_n^{(0)}}^{(2)}(t)$.
The functions $u_n(x,t)$ satisfy the initial conditions
\begin{equation}
\label{eq36}
u_n (x,t)\vert_{t=0}=\varphi^{(0)}_n(x)=v_n^{(0)}(x,0,{\bf C}_n),
\end{equation}
where $v_n^{(0)}(x,0,{\bf C}_n)$ is given by (\ref{eq33}).

The general solution of the EE system (\ref{sys-111})--(\ref{sys-333}) can be presented in the Cauchy form; that is,
the initial conditions ${\bf C}=(\sigma (0, D), x (0, D), \alpha^{(0)}(0,D))$ can be taken as integration constants (see, e.g., (\ref{sys-444}) and (\ref{sys-4444})).
Then the constants ${\bf C}_n$ in (\ref{eq36}) are determined from the equalities
\begin{equation}
{\bf C}=\big(\sigma(0,D),x(0,D),\alpha^{(2,2)}(0,D)\big)={\bf C}_n=\big(\sigma_n,x_n,\alpha_n^{(2)}\big),\label{44aq1q}
\end{equation}
where
\begin{eqnarray}
&&\sigma_n=\int\limits^{+\infty}_{-\infty}v_n^{(0)}(x,0)dx; \label{44aq1qA}\\
&&x_{n}=\frac1{\sigma_{n}}\int\limits^{+\infty}_{-\infty}xv_{n}^{(0)}(x,0)dx;\label{44aq1q1}\\
&&\alpha^{(2)}_{n}=\frac1{\sigma_{n}}\int\limits^{+\infty}_{-\infty}\Delta x^2v_{n}^{(0)}(x,0)dx.\label{44aq1q2}
\end{eqnarray}

Let us calculate the constants ${\bf C}_n$. To construct asymptotic solutions accurate to $O(D^{3/2})$, it suffices to consider the case $M= 2$:
\begin{align}
\sigma_n&=\int\limits^{+\infty}_{-\infty}v_n^{(0)}(x,0)dx= \int\limits^{+\infty}_{-\infty}\frac{(-1)^n}{\sqrt{2^nn!}}
\Big[\frac{Z^{(+)}(t)}{Z^{(-)}(t)}\Big]^{n/2}H_n
\left(\sqrt{\frac{b}{DZ^{(-)}(t)Z^{(+)}(t)}}\Delta x\right)\times\nonumber\\
&\times N_D\sqrt{-\frac{W^{(-)}(t)}{bZ^{(-)}(t)}}
\frac{\sigma^{(M)}(t,{\bf C})}{\sigma(0)}\exp\Big[
\frac{1}{2D}\frac{W^{(-)}(t)}{Z^{(-)}(t)}\Delta x^2\Big]\bigg|_{t=0}dx=\nonumber\\
&=\frac{(-1)^nN_D}{\sqrt{2^nn!}}\int\limits^{+\infty}_{-\infty}
H_n\left(\sqrt{\frac{b}{D}}\Delta x\right)
\exp\Big[-\frac{b}{2D}\Delta x^2\Big]dx=\nonumber\\
&=\frac{(-1)^nN_D}{\sqrt{2^nn!}}\sqrt{\frac{D}{b}}\int\limits^{+\infty}_{-\infty}
H_n\left(y\right)
\exp\Big[-\frac{1}{2}y^2\Big]dy=\nonumber\\
&=\sqrt[4]{\frac{D}{\pi b}}\begin{cases}0,&n=2l+1;\\
\dfrac{\sqrt{2\pi}}{\sqrt{2^{2l}(2l)!}}
\dfrac{(2l)!}{l!},&n=2l.\end{cases}\label{44aq}
\end{align}
Hence, $\sigma_{2l+1}\equiv0$ and $\sigma_{2l}=\sqrt[4]{\frac{4\pi D}{ b}}\frac{\sqrt{(2l)!}}{2^{l}l!}$.
Here we  use the relation
\begin{equation}
H_n(y)=\frac{\partial^n}{\partial s^n}\exp\Big[-s^2+2sy \Big]\bigg|_{s=0},\label{44aq1q2-1}
\end{equation}
which follows from the definition of the Hermite polynomials \cite{Beitman2}:
\[\exp\Big[-s^2+2sy \Big]=\sum_{n=0}^\infty\frac{s^{n}}{n!}H_n(y).\]
In view of the relation \eqref{44aq1q2-1} we write
\begin{align}
I_n&=\int\limits^{+\infty}_{-\infty}
H_n(y)
\exp\Big[-\frac{1}{2}y^2\Big]dy=\nonumber\\
&=\frac{\partial^n}{\partial s^n}\int\limits^{+\infty}_{-\infty}\exp\Big[-\frac{1}{2}y^2\Big]
\exp\big[-s^2+2sy\big]dy\bigg|_{s=0}=\nonumber\\
&=\frac{\partial^n}{\partial s^n}
\exp\big[-s^2\big]\sqrt{2\pi}\exp\big[2s^2\big]\bigg|_{s=0}
=\sqrt{2\pi}\frac{\partial^n}{\partial s^n}
\sum_{k=0}^\infty\frac{s^{2k}}{k!}\bigg|_{s=0}=\nonumber\\
& =\sqrt{2\pi}\sum_{k=n}^\infty\frac{2k(2k-1)\cdots(2k-n+1)}{k!}s^{2k-n}\bigg|_{s=0}
=\begin{cases}0,&n=2l+1;\\
\sqrt{2\pi}\dfrac{(2l)!}{l!},&n=2l,\end{cases}\label{44aq2}
\end{align}
and
\begin{equation}
\int\limits^{+\infty}_{-\infty}
\exp\Big[-\varpi^2y^2+2sy\Big]dy=\frac{\sqrt{\pi}}{\varpi}\exp\Big[\frac {s^2}{\varpi^2}\Big],\quad N_D=\sqrt[4]{\frac{b}{\pi D}}.\label{44aq1}
\end{equation}

Similarly, we obtain
\begin{align}
x_{2n}&=\frac1{\sigma_{2n}}\int\limits^{+\infty}_{-\infty}xv_{2n}^{(0)}(x,0)dx=\nonumber\\
&=\frac1{\sigma_{2n}}\int\limits^{+\infty}_{-\infty}\Delta xv_{2n}^{(0)}(x,0)dx+\frac{x_0}{\sigma_{2n}}\int\limits^{+\infty}_{-\infty}v_{2n}^{(0)}(x,0)dx=x_0\label{44aq3}
\end{align}
and
\begin{align}
\alpha^{(2)}_{2n}&=\frac1{\sigma_{2n}}\int\limits^{+\infty}_{-\infty}\Delta x^2v_{2n}^{(0)}(x,0)dx=\nonumber\\
&=\frac1{\sigma_{2n}}\frac{N_D}{\sqrt{2^{2n}(2n)!}}\int\limits^{+\infty}_{-\infty}\Delta x^2
H_{2n}\left(\sqrt{\frac{b}{D}}\Delta x\right)
\exp\Big[-\frac{b}{2D}\Delta x^2\Big]dx=\nonumber\\
&=\frac1{\sigma_{2n}}\frac{N_D}{2^n\sqrt{(2n)!}}\sqrt{\frac{D^3}{b^3}}\int\limits^{+\infty}_{-\infty}y^2
H_{2n}(y)
\exp\Big[-\frac{1}{2}y^2\Big]dy=\nonumber\\
&= \sqrt[4]{\frac{ b}{4\pi D}}\frac{2^{n}n!}{\sqrt{(2n)!}}\sqrt[4]{\frac{b}{\pi D}}\frac{1}{2^n\sqrt{(2n)!}}\sqrt{\frac{D^3}{b^3}}\sqrt{2\pi}\frac{(2n)!}{n!}
[1+4n]= 
\frac{D}{b}[1+4n].\label{44aq6}
\end{align}
Here we take into account that  $H'_{n}(y)=2nH_{n-1}(y)$  and
\[
J_n=\int\limits^{+\infty}_{-\infty}y^2
H_{2n}(y)
\exp\Big[-\frac{1}{2}y^2\Big]dy.\]
Integrating this by part and setting $U=yH_{2n}(y)$, $dU=\bigl(yH_{2n}(y)\bigr)'dy$, $dV=ye^{-y^2/2}dy$ and
$V=-e^{-y^2/2}$, we obtain
\begin{eqnarray}
J_n&=&-yH_{2n}(y)e^{-y^2/2}\bigg|^\infty_{-\infty}+\int\limits^{+\infty}_{-\infty}\big[y
H'_{2n}(y)+H_{2n}(y)\big]
\exp\Big[-\frac{1}{2}y^2\Big]dy=\nonumber\\
& =&I_{2n}+\int\limits^{+\infty}_{-\infty}y
(4n)H_{2n-1}(y)
\exp\Big[-\frac{1}{2}y^2\Big]dy.\nonumber
\end{eqnarray}
Integrating this by part and setting $U=H_{2n-1}(y)$, $dU=H'_{2n-1}(y)dy$, $dV=ye^{-y^2/2}dy$ and
$V=-e^{-y^2/2}$, we obtain
\begin{eqnarray}
J_n&=&I_{2n}+4n\int\limits^{+\infty}_{-\infty}
H'_{2n-1}(y)\exp\Big[-\frac{1}{2}y^2\Big]dy=I_{2n}
+4(2n)(2n-1)I_{2n-2}=\nonumber\\
&=&\sqrt{2\pi}\Big[\frac{(2n)!}{n!}+4(2n)(2n-1)\frac{(2n-2)!}{(n-1)!}\Big]
=\sqrt{2\pi}\frac{(2n)!}{n!}[1+4n].\label{44aq5}
\end{eqnarray}

Thus, under the conditions
\begin{equation}
\sigma_{2n}(0,D)=\sqrt[4]{\frac{4\pi D}{ b}}\frac{\sqrt{(2n)!}}{2^{n}n!},\; x_{2n}(0,D)=x_0,\;
\alpha^{(2)}_{2n}(0,D)=\frac{D(1+4n)}{b},\label{44aq1A}
\end{equation}
the functions $u_{2n+1}(x,t,D)\equiv 0$, and $u_{2n}(x,t,D)$ determined by equations  (\ref{solution-family})
are the desired, accurate to $O(D^{3/2})$, asymptotic solutions for the nonlinear equation (\ref{eq0}).

\section{Large time asymptotics to quasi-stationary solutions}
\label{quasistat}

We have already noted (see Section \ref{asslinprob}) that the semiclassical asymptotics constructed above (see equations (\ref{eq15})) hold for a finite time interval $t\in [0, T]$, where $T$ does not depend on $D$ and the asymptotic estimate fails to hold for large $T$.

For certain types of the functions $a(x, t)$, $b(x, y, t)$, $V(x, t)$, and $W(x,y,t)$ entering into equation (\ref{eq0}), the residual of asymptotic solutions may be limited or even decrease with time \cite{12c}. For example, from (\ref{sys-4444}) it follows  that the variance increases  infinitely as $T\to\infty$, and therefore $u(x,t)$ tends to zero. As $u(x,t)=0$ is an exact stationary solution of (\ref{eq0}), the constructed semiclassical asymptotics   (\ref{eq15}) can be treated as a perturbation of the stationary solution. This raises the problem of finding other exact solutions with the asymptotic solutions treated as their perturbations.

Consider an example of the one-dimensional Fisher--KPP equation (\ref{eq0}) with the following coefficients:
\begin{eqnarray}
& V(x,t)=0, \quad W(x,y,t)=0,\quad x\in {\mathbb R}, \quad a(x,t)=a={\rm const},\quad a>0,\nonumber\\
& b(x,y,t)=b(x-y), \quad \displaystyle\int\limits_{-\infty}^{\infty}b(x-y)dy=B={\rm const}, \quad |B|<\infty . \label{anal1}
\end{eqnarray}

Equation (\ref{eq0}) in this case takes the form
\begin{equation}
u_t(x,t)=Du_{xx}(x,t)+au(x,t)-\varkappa u(x,t)\int_{-\infty}^{\infty} b(x-y)u(y,t)dy. \label{anal2}
\end{equation}

A spatially homogeneous particular solution of equation (\ref{anal2}) is sought in the form
\begin{equation}
u(x,t) = u_0(x,t) =\beta(t),\quad u_0(x,t)\big|_{t=0}= \beta\big|_{t=0}=\beta_{0}. \label{anal_resh4z}
\end{equation}
Substituting (\ref{anal_resh4z}) in (\ref{anal2}), we obtain the Verhulst equation for the function $\beta (t)$:
\begin{equation}
\dot \beta =a \beta-{\varkappa} B \beta^2. \label{anal10q}
\end{equation}
Using the procedure by which we have obtained
(\ref{sys-112}) and  (\ref{sys-4444}), we find
\begin{equation}
u_0(x,t) =\beta(t) = \frac{\beta_{0}e^{at}}{1+\varkappa\frac{\beta_{0}}a B( e^{at}-1)},  \label{anal_resh4aa}
\end{equation}
whence
\begin{equation}
\lim_{t\to\infty} u_0(x,t) = \frac{a}{\varkappa B}=u_\infty.   \label{anal_resh14}
\end{equation}

We seek a solution to equation (\ref{anal2}) as
\begin{equation}
u(x,t) = u_0(x,t) + \bar u(x,t).    \label{anal_resh14a}
\end{equation}
We refer to solutions of this type as quasi-stationary functions. Note that  solution (\ref{anal_resh4z}) does not belong to $L_2$,
and therefore it can not be described by semiclassical asymptotics (\ref{eq15}).

Substituting (\ref{anal_resh14a}) in (\ref{anal2}), we obtain
\begin{gather}
\bar  u_t(x,t)-D\bar u_{xx}(x,t)-\big[a-\varkappa B\beta(t)\big]\bar u (x,t) + \nonumber\\
+ \varkappa \beta(t)\int\limits_{-\infty}^{\infty} b(x-y)\bar u(y,t)dy
+\varkappa\bar u(x,t)\int\limits_{-\infty}^{\infty} b(x-y)\bar u(y,t)dy=0.\label{anal_resh18bbb1}
\end{gather}

Equation (\ref{anal_resh18bbb1}) has the structure of equation (\ref{eq0}) and can be solved by the above semiclassical asymptotic method.

Consider the problem of small perturbations for the solutions (\ref{anal_resh4aa}). Let the initial condition for
the function $\bar u(x,t)$ be
\begin{equation}
\bar u(x,0)=\varepsilon \varphi(x),\quad 0<\varepsilon\ll\frac{\beta_{0}}{\|\varphi(x)\|}\,.
\label{vst_beta0}
\end{equation}

Note that the small parameter $\varepsilon $ in (\ref{vst_beta0}) can be chosen in various ways.

For instance, if we put  $t=\varepsilon s$,  where $s$ is  dimensionless time and $T$ is
the characteristic evolution time of a population system, we can choose $\varepsilon=1/T$ and the asymptotic solutions
with the small parameter $\varepsilon$ will describe the behavior of the system
for a large time $T$.

Here the function $\varphi(x)$ belongs to the Schwartz space.
For equation (\ref{anal2}),
the Cauchy problem (\ref{anal_resh18bbb1}), (\ref{vst_beta0}) generates the initial condition
\begin{equation}
u\Big|_{t=0}=\beta_{0} + \varepsilon\varphi(x).
\label{vst_beta00}
\end{equation}

We  seek an asymptotic solution to equation (\ref{anal_resh18bbb}) in the class of functions regularly depending on the small parameter $\varepsilon$, $\varepsilon \to 0$:
\begin{equation}
\bar u(x,t) = \bar u (x, t,\varepsilon) = \varepsilon\bar u^{(1)} (x,t) +
\varepsilon^2\bar u^{(2)} (x, t)+\ldots\, .
\label{anal_resh15}
\end{equation}

Substituting the expansion (\ref{anal_resh15}) in (\ref{anal2}) and equating terms of equal powers of $\varepsilon$, we obtain
\begin{eqnarray}
&&
\bar  u_t^{(1)}(x,t)-D\bar u_{xx}^{(1)}(x,t)-\big[a-\varkappa B\beta(t)\big]\bar u^{(1)} (x,t)+ \nonumber\\
&&\quad+
\varkappa \beta(t)\int\limits_{-\infty}^{\infty} b(x-y)\bar u^{(1)}(y,t)dy=0, \label{anal_resh18a} \\
&&\bar  u_t^{(2)}(x,t)-D\bar u_{xx}^{(2)}(x,t)+
\varkappa \beta(t)\int\limits_{-\infty}^{\infty} b(x-y)\bar u^{(2)}(y,t)dy-\nonumber\\
&&\quad-\big[a-\varkappa B\beta(t)\big]\bar u ^{(2)}(x,t) +\varkappa\bar u^{(1)}(x,t)\int\limits_{-\infty}^{\infty} b(x-y)\bar u^{(1)}(y,t)dy=0,\label{anal_resh18b} \\
&&\dots\dots\dots\dots\dots\dots\dots \nonumber
\end{eqnarray}

The initial conditions (\ref{vst_beta0}) yield
\begin{equation}
\bar u^{(1)}\Big|_{t=0}=\varphi(x), \quad \bar u^{(k)}\Big|_{t=0}=0, \quad  k=2,3,\ldots .
\label{vst_beta0q}
\end{equation}

Applying the Fourier transform
$$
\tilde f(\omega)=\frac{1}{\sqrt{2\pi}}\int\limits_{-\infty}^{\infty}f(x)e^{-i\omega x} d x, \quad
 f(x)=\frac{1}{\sqrt{2\pi}}\int\limits_{-\infty}^{\infty}\tilde f(\omega)e^{i\omega x} d \omega
$$
to the right and left sides of equation (\ref{anal_resh18a}), in view of the initial conditions (\ref{vst_beta0q}), we find
\begin{eqnarray}
&& \tilde u^{(1)}_t(p,t) +\left\{D p^2 - a+\varkappa \big[B +\sqrt{2\pi}\,\, \tilde b(p)\big] \beta(t) \right\}\tilde u^{(1)}(p,t)=0,\label{anal_resh18bbb}\\
&&\tilde u^{(1)}(p,t)\Big|_{t=0}=\tilde\varphi(p)\nonumber .
\end{eqnarray}

Here, $\tilde u^{(1)}$, $\tilde b$, and $\tilde \varphi$ are the  Fourier images of the functions $\bar u^{(1)}$, $b$, and $\varphi$, respectively. From (\ref{anal_resh18bbb}) we find
\begin{equation}
\tilde u^{(1)}(p,t)=\tilde\varphi(p)\exp\left\{-D p^2 t + at -\varkappa\big[B+ \sqrt{2\pi}\,\,\tilde b(p)\big]\chi(t)\right\},\label{anal_resh18bbbb}
\end{equation}
where $\chi(t)=\int_0^t \beta(s)ds$.

Applying the inverse Fourier transform to $\tilde u^{(1)}(p,t)$ and taking into account (\ref{vst_beta0q}), we write
\begin{equation}
\bar u^{(1)}(x,t)= \int\limits_{-\infty}^\infty  G^{(1)}(x,y,t,0) \varphi(y) dy, \label{anal_resh18aq}
\end{equation}
where
\begin{align}
G^{(1)}(x,y,t,s)&=\frac1{2\pi}\int\limits_{-\infty}^\infty\exp\Big\{ i p(x-y) + a(t-s) - Dp^2 (t-s)- \nonumber\\
&-\varkappa \big[B+\sqrt{2\pi}\, \tilde b(p)\big]\,[\chi (t)-\chi (s)] \Big\}dp. \label{anal_resh18bq}
\end{align}

Similarly, from (\ref{anal_resh18b}), (\ref{vst_beta0q}), and (\ref{anal_resh18bq}), we obtain
\begin{equation}
\bar u^{(2)}(x,t) = -\int\limits_{0}^{t} d s
\int\limits_{-\infty}^{\infty} G^{(1)}(x,y,t, s) f(y,s)\, d y.
\label{anal_resh20}
\end{equation}
From (\ref{anal_resh18b}) it follows that
\begin{equation}
f(y,t) = \varkappa\bar u^{(1)}(y,t)\int\limits_{-\infty}^{\infty}
b(y-z)\bar u^{(1)}(z,t)\, dz .
\end{equation}

In the same way, we can find the functions $\bar u^{(k)}$, $k={3, 4,\dots  }$.
Substituting (\ref{anal_resh18bq}) in (\ref{anal_resh18aq}), we have
\begin{align}
\bar u^{(1)}(x,t)&= \int\limits_{-\infty}^\infty  \frac{1}{2\pi}
\int\limits_{-\infty}^{\infty} d p
\exp\Big[i p (x-y) -Dp^2t+at- \nonumber\\
&-\varkappa \big(B+\sqrt{2\pi}\,\,\tilde b(p)\big)\chi(t)\Big]  \varphi(y) dy=\nonumber\\
&=\frac{1}{\sqrt{2\pi}}
\int\limits_{-\infty}^{\infty} d p
\exp\Big\{i px -Dp^2t+at- 
\varkappa \big[B+\sqrt{2\pi}\,\tilde b(p)\big]\chi(t)\Big\} \tilde \varphi(p).\label{anal_resh18aq1}
\end{align}
As an example, consider an influence function and initial functions of the form
\begin{equation}
b(x)=b_0\exp\big[-x^2/\gamma^2\big]
\label{anal_resh18aq2}
\end{equation}
and
\begin{equation}
\varphi(x)= \varphi_n(x)=N u_n\big(\theta(x-x_0)\big),
\label{anal_resh18aq3}
\end{equation}
respectively.
Here $u_n(x)$ are the Hermite functions, $n=0,1, \dots$,  $(\theta  >0$) and $x_0$ are parameters.
Then we have
$$
\tilde b(p)=\frac{\gamma}{\sqrt2} b_0\exp\Big[-\frac {p^2\gamma^2}4\Big], \quad
\tilde \varphi_n(p)=(-i)^n\frac N\theta e^{-ipx_0}u_n\Big(\frac p\theta\Big).
$$
Substituting (\ref{anal_resh18aq2}) and (\ref{anal_resh18aq3}) in (\ref{anal_resh18aq1}), we obtain
\begin{eqnarray}
&&\bar u^{(1)}(x,t)=\bar u^{(1)}_n(x,t)=\nonumber\\
&&\quad=\frac{1}{\sqrt{2\pi}}(-i)^n\frac N\theta
\int\limits_{-\infty}^{\infty} d p\,
\exp\Big[i p(x-x_0) -\big(Dp^2-a\big)t- \varkappa B\chi(t)\Big]\times\nonumber\\
&&\quad\times\sum_{k=0}^\infty\frac1{k!}\Big[\varkappa b_0\sqrt{\pi}\gamma \chi(t)\Big]^k\exp\Big[-k\frac{p^2\gamma^2}4\Big] u_n\Big(\frac p\theta\Big)=\nonumber\\
&&\quad=\frac{1}{\sqrt{2\pi}}(-i)^n\frac N\theta
\exp\big[ at- \varkappa B\chi(t)\big]\sum_{k=0}^\infty\frac1{k!}\Big[\varkappa b_0\sqrt{\pi}\gamma \chi(t)\Big]^k\times\nonumber\\
&&\quad\times\int\limits_{-\infty}^{\infty} d p\,\exp\Big[i p(x-x_0)-\Big(Dt+k\frac{\gamma^2}4+\frac1{2\theta^2}\Big)p^2\Big]\frac1{\sqrt{2^nn!\sqrt\pi}}  H_n\Big(\frac p\theta\Big)=\nonumber\\
&&\quad=(-i)^n\frac{N}{\sqrt{2\pi}}
\exp\Big[ at- \varkappa B\chi(t)\Big]\sum_{k=0}^\infty\frac1{k!}\Big[\varkappa b_0\sqrt{\pi}\gamma \chi(t)\Big]^k\times\nonumber\\
&&\quad\times\Big(\frac{4\pi}{4D\theta^2t+k\gamma^2\theta^2+2}\Big)^{1/2}
\Big(\frac{4D\theta^2t+k\gamma^2\theta^2-2}{4D\theta^2t+k\gamma^2\theta^2+2}\Big)^n\times\nonumber\\
&&\quad\times\frac{1}{\sqrt{2^nn!\sqrt\pi}} \exp\Big[\frac{-\theta^2(x-x_0)^2}{4D\theta^2t+k\gamma^2\theta^2+2}\Big] H_n\Big(\frac {2i\theta(x-x_0)}{[(4Dt+k\gamma^2)^2\theta^4-4]^{1/2}}
\Big).
\label{anal_resh18aq4}
\end{eqnarray}
Expand the functions $\bar u^{(1)}(x,t)$   from  (\ref{anal_resh18aq4}) in series of the Hermite functions $u_m\big(\theta(x-x_0)\big)$:
\begin{equation}
\bar u^{(1)}_n(x,t)= \sum_{m=0}^\infty C_m^{(n)}(t)u_m\big(\theta(x-x_0)\big). \label{anal_resh18aq41}
\end{equation}
Using
$$
\int\limits_{-\infty}^{\infty} d x\, u_m\big(\theta(x-x_0)\big)\exp\Big[i p(x-x_0)\Big]=\sqrt{2\pi}\frac {(+i)^m}\theta u_m\Big(\frac p\theta\Big),
$$
we have
\begin{align}
C_m^{(n)}(t)&=(-i)^{n-m}\frac N{\theta^2}\int\limits_{-\infty}^{\infty} d p\,
\exp\Big[\big(-Dp^2+a\big)t- \varkappa B\chi(t)\Big]\times\nonumber\\
&\times\sum_{k=0}^\infty\frac1{k!}\Big[\varkappa b_0\sqrt{\pi}\gamma \chi(t)\Big]^k\exp\big[-k\frac{p^2\gamma^2}4\big] u_n\Big(\frac p\theta\Big)u_m\Big(\frac p\theta\Big)=\nonumber\\
&=(-i)^{n-m}\frac N{\theta^2}
\exp\Big[ at- \varkappa B\chi(t)\Big]\sum_{k=0}^\infty\frac1{k!}\Big[\varkappa b_0\sqrt{\pi}\gamma \chi(t)\Big]^k
\frac \theta{\sqrt{2^{n+m}n!m!\pi}}\times\nonumber\\
&\times \int\limits_{-\infty}^{\infty} d s\,\exp\big[-\big(D\theta^2t+k\frac{\gamma^2\theta^2}4+1\big)s^2\big] H_n(s)H_m(s).
\label{anal_resh18aq45}
\end{align}

For the particular case $n=0$, from (\ref{anal_resh18aq45}) we obtain  $C_{2l+1}^{(0)}(t)=0$ and
\begin{align}
C_{2l}^{(0)}(t)&=(-1)^{l}\frac N{\theta}
\exp\Big[at- \varkappa B\chi(t)\Big]\sum_{k=0}^\infty\frac1{k!}\Big[\varkappa b_0\sqrt{\pi}\gamma \chi(t)\Big]^k\times\nonumber\\
&\times\frac 1{2^l\sqrt{(2l)!\pi}} \int\limits_{-\infty}^{\infty} d s\exp\big[-\big(D\theta^2t+k
\frac{\gamma^2\theta^2}4+1\big)s^2\big] H_{2l}(s)=\nonumber\\
&=\frac{(-1)^{l} N}{\theta 2^l\sqrt{(2l)!\pi}}
\exp\Big[ at- \varkappa B\chi(t)\Big]\sum_{k=0}^\infty\frac1{k!}\Big[\varkappa b_0\sqrt{\pi}\gamma \chi(t)\Big]^k\times\nonumber\\
&\times \frac{\sqrt{\pi}\big(D\theta^2t+k\frac{\gamma^2\theta^2}4\big)^l}{\big(D\theta^2t+
k\frac{\gamma^2\theta^2}4+1\big)^{(2l+1)/2}}H_{2l}(0)=\nonumber\\
&=\frac{ N\sqrt{(2l)!}}{\theta 2^{l-1}l! }
\exp\big[at- \varkappa B\chi(t)\big]\sum_{k=0}^\infty\frac1{k!}\Big[\varkappa b_0\sqrt{\pi}\gamma \chi(t)\Big]^k\times\nonumber\\
&\times \frac1{\sqrt{4D\theta^2t+k\gamma^2\theta^2+4}}\bigg(\frac{4D\theta^2t+k\gamma^2\theta^2}
{4D\theta^2t+k\gamma^2\theta^2+4}\bigg)^{l}.
\label{anal_resh18aq42}
\end{align}

Direct verification shows that we have $C_{0}^{(0)}(0)=1$ for $t\to 0$, $C_{2l}^{(0)}(0)=0$ for $l > 0$, and $C_{2l}^{(0)}(t)\neq 0$ for $t>0$.
Thus, even if the initial  distribution $ \bar u(x,0)$ in (\ref{vst_beta0}) is unimodal, the solution  $\bar u(x,t)$  becomes multimodal for large $t$, as follows from (\ref{anal_resh18aq41})
and (\ref{anal_resh18aq42}).

In particular, for $t\to\infty$ we have
\begin{equation}
C_{2l}^{(0)}(t)\thicksim \frac{ N\sqrt{(2l)!}}{ 2^{l}l! \theta^2\sqrt{Dt}}
\exp\Big[\frac{b_0\pi at}{\gamma B} \Big],\quad t\to\infty. \label{anal_resh18aq6}
\end{equation}
Here, $\chi(t)\thicksim \displaystyle\frac{at}{\varkappa B}$ as $t\to\infty$.

 From the properties (\ref{anal_resh18aq45})--(\ref{anal_resh18aq6})
of the expansion (\ref{anal_resh18aq41}),  we can conclude that the  asymptotic solutions (\ref{anal_resh18aq4})
describe the evolution of spatial patterns. Thus, the asymptotic approach used
can be considered a contribution to the analytical methods of studying pattern formation with the use of the Fisher--KPP equation.

\section{Conclusion remarks }
\label{concusion}

 In this paper, we propose a modified Maslov complex germ method to construct real solutions to the 1D non-local Fisher--KPP equation (\ref{eq0}).

Here, the semiclassical approximation is considered in the class   $\mathcal{P}_t^D $ of  trajectory-concentrated functions (\ref{eq1}).
The key point of the formalism developed   is the  dynamic Einstein--Ehrenfest system describing the evolution
 of  the moments.
The general solution of the EE system, on the one hand,  allows one to associate a linear equation with the original Fisher--KPP equation.
On the other hand, we formulate algebraic conditions (\ref{algeb-cond-1}) for the integration constants.
As a result, we obtain a set of functionals which are integrals of the Fisher--KPP equation.

The original nonlocal Fisher--KPP equation is solved by substituting the found integrals in the solution
of the associated linear equation. In this paper, we focus on the most difficult part of the asymptotic
series, namely the search for leading terms.

The solution of the linear associated equation (\ref{eq23a}) is found by using a set of symmetry operators
   $\big(\hat{a}^{(\pm)}(t) \big)$  of the form (\ref{eq30}), (\ref{eq30aa}).
The higher-order terms are determined by solving
a linear problem, and this does not cause principal difficulties.

The family of particular asymptotic solutions (leading terms) constructed in explicit form in Section \ref{leadsemiclsol}
can be considered to possess a completeness property in the following sense: At the time zero, $t=0$, these solutions form a basis in $L_2$. Using this basis, we can find  solutions
of the associated linear equation for a wide class of initial distributions. Then the solution
of the nonlinear Fisher--KPP equation will be reduced to algebraic conditions for the  integration constants
of the dynamic EE system.

However, the semiclassical asymptotics are valid for a finite time interval $[0,T]$.
Estimation of this interval is  beyond the scope of the semiclassical method and requires additional study.
Note that in particular cases
this interval can be quite small, i.e. the asymptotics will be valid at small times.
The smallness of $T$ implies that a spatial structure (pattern) has no time to form in the interval $[0,T]$.

To go beyond these limitations,  in the final Section \ref{quasistat},
we construct asymptotics  representing small perturbations on the background of
an  exact quasi-stationary solution of the nonlocal Fisher--KPP equation. They are different from
the semiclassical asymptotics and can describe the evolution of  distributions for  large times ($T\to\infty$).
In the example considered, an initial  ($t=0$) unimodal  distribution becomes mutimodal at large $t$, and this can be treated as pattern formation.
In other words, the developed approach can be used to analyze the dynamics of pattern formation.

The semiclassical approximation considered here for the 1D nonlocal Fisher--KPP equation
in the class of trajectory concentrated functions  $\mathcal{P}_t^D$
provides  the basic formalism for further generalization to multidimensional cases.

Multidimensional asymptotics can be sought for function spaces different from $\mathcal{P}_t^D$, and this offers additional ways to study the pattern formation in the population dynamics.
For instance, we used a semiclassical formalism to describe the formation of 2D patterns by reducing an original 2D problem to a 1D problem for a one-dimensional curve in a two-dimensional space, no invoking exact solutions of the 2D Fisher--KPP equation \cite{jphys2014}.
However, the solution of the resulting 1D equations gives a complete picture of the population dynamics. Therefore, in our opinion, further generalization of the proposed methods is a topical problem.

\section*{Acknowledgements}
The work was supported in part by Tomsk State University under the International Competitiveness
Improvement Program and by Tomsk Polytechnic University under the International Competitiveness
Improvement Program.

\end{document}